\documentstyle[11pt]{article}

\newtheorem{th}{Theorem}[section]

\newtheorem{ex}[th]{Example}
\newtheorem{Theorem}[th]{Theorem}
\newtheorem{Lemma}[th]{Lemma}
\newtheorem{Proposition}[th]{Proposition}
\newtheorem{Remark}[th]{Remark}
\newtheorem{Remarks}[th]{Remarks}
\newtheorem{Corollary}[th]{Corollary}

\newtheorem{ccote}[th]{}

\newcommand{\decale}[1]{\par\noindent\hskip 3em\llap{#1\enspace}\ignorespaces}

\newcommand{\SS}{\S \kern .2em}
\newcommand{\hfb}{\newline}

\newcommand{\preu}{{\sc Proof: \ }}
\newcommand{\cit}[1]{\ref{#1}}

\newcommand{\fl}[1]{\buildrel{#1}\over{\longrightarrow}}

\newcommand{\cqfd}{\unskip\kern 6pt\penalty 500
\raise -2pt\hbox{\vrule\vbox to10pt{\hrule width
 4pt\vfill\hrule}\vrule}\smallskip}

\newcommand{\bbr}{{\bf R}}
\newcommand{\bbc}{{\bf C}}
\newcommand{\bbz}{{\bf Z}}
\newcommand{\bbq}{{\bf Q}}
\newcommand{\bbn}{{\bf N}}

\newcommand{\calb}{{\cal B}}

\newcommand{\calf}{{\cal F}}

\newcommand{\cali}{{\cal I}}
\newcommand{\calj}{{\cal J}}

\newcommand{\call}{{\cal L}}

\newcommand{\calr}{{\cal R}}
\newcommand{\cals}{{\cal S}}

\newcommand{\calu}{{\cal U}}
\newcommand{\calv}{{\cal V}}

\newcommand{\pol}{{\rm Pol\,}}
\newcommand{\apol}{{\rm APol\,}}
\newcommand{\upol}{{\rm UP\,}}
\newcommand{\vpol}{{\rm VP\,}}
\newcommand{\polr}{{\rm Pol_\bbr\,}}
\newcommand{\apolr}{{\rm APol_\bbr\,}}
\newcommand{\upolr}{{\rm UP_\bbr\,}}

\newcommand{\rep}{\bbr_{+}}

\newcommand{\union}{\cup}
\newcommand{\pcirc}{\kern .7pt {\scriptstyle \circ} \kern 1pt}

\newcommand{\plaint}[1]{}
\newcommand{\iso}{\cong}
\newcommand{\Ann}{{\rm Ann}}

\newcommand{\btab}{\begin{tabbing}}
\newcommand{\etab}{\end{tabbing}}

\title{The cohomology ring of polygon spaces}
\author{Jean-Claude HAUSMANN
\and Allen KNUTSON
\footnote{Both authors thank the Fonds National Suisse de la Recherche Scientifique for its support.}  }
\date{}

\begin{document}    \maketitle

\begin{abstract}

We compute the integer cohomology rings of the ``polygon spaces'' introduced
in \cite{Ha,Kl,KM}. This is done by embedding them in certain toric varieties;
the restriction map on cohomology is surjective and we calculate its kernel 
using ideas from the theory of Gr\"obner bases.
Since we do not invert the prime 2, we can tensor with $\bbz_2$;
halving all degrees we show this produces 
the $\bbz_2$ cohomology rings of the planar polygon spaces.
In the equilateral case, where there is an action of the symmetric group
permuting the edges, we show that the induced action on the integer cohomology 
is {\em not} the standard one, despite it being so on the rational cohomology 
\cite{Kl}. 
Finally, our formulae for the Poincar\'e polynomials are more 
computationally effective than those known \cite{Kl}.
\end{abstract}

\section*{Introduction} \label{intro}

A ``polygon space'' $\pol(\alpha_1,\alpha_2,...\alpha_m),
\alpha_i \in {\bbr}_+$
can be seen to arise in several ways:

1. the family of piecewise linear paths in ${\bbr}^3$, whose $i$th step (which
is of length $\alpha _i$) can proceed in any direction subject to the polygon
ending where it begins, considered up to rotation and translation

\newcommand{\Sym}{{\rm Sym}}
2. the ``semistable'' configurations of $m$ weighted points in ${\bbc P}^1$,
the $i$th of weight $\alpha _i$, where a configuration is considered unstable
if more than half the total weight sits in one place, modulo M\"obius
transformations. This description is of the most classical interest \cite{DM},
particularly in the equal-weight case \cite{Popp},
since the quotient by $\Sym_m$
is a compactification of the moduli space of $m$-pointed
genus zero curves

3. (when the $\alpha _i$ are integral) the geometric invariant theory quotient
of the Grassmannian of $2$-planes in ${\bbc}^n$ by $T^n$, where the $\{\alpha
_i\}$ specify an action of $T^n$ on the canonical bundle.

The connection of the first to the second is made in \cite{Kl} and \cite{KM};
the second to the third in \cite{GM}, \cite{GGMS} and the first to the third in
\cite{HK}. This paper draws much from the polygonal
 intuition and will
concentrate on the first.

In this paper we compute the integer cohomology rings of these spaces, in the
(generic) case that they are smooth. There are partial results in the
literature. Klyachko \cite{Kl} showed that the cohomology groups were
torsion-free and calculated their rank. Brion \cite{Br}
calculated the rational cohomology ring in the equal-weight case with an
odd number of sides modulo the symmetric group, Kirwan \cite{KiPoly}
the de Rham cohomology ring in the ordered equal-weight case.
It seems that a slight refinement of Kirwan's method would
require that one only invert the prime 2.

Our approach is very different, and makes heavy use of toric varieties, whose
integer cohomology rings are known by the theorem of Danilov. While a polygon
space is not (usually) a toric variety itself, it embeds in one in a very
special way: as a transverse self-intersection of a toric subvariety.

This gives a map from the cohomology ring of the ambient toric variety,
our {\em upper path space}, to that of the polygon space itself.
We then have four tasks to complete:

1. Compute the cohomology ring of the upper path space,
using (a mild extension of) Danilov's result.

2. Show that the restriction map on cohomology is surjective.

3. Show that the kernel of this map is the annihilator of the
Poincar\'e dual of the submanifold.

4. Compute the annihilator.

The first is a very polygon-theoretic argument, and is in Section \ref{45}. We
prove the second and third as part of a more general study of {\em
even-cohomology spaces}, based on the  fact that
 $H^{\rm odd} = 0$ for not only the
polygon space and upper path space, but also the difference.\footnote{%
In fact claims 2 and 3 hold much more generally, and are the basis
of Shaun Martin's unpublished but very influential $G$-to-$T$ argument
\cite{SKMthesis}.
Owing to our very special even-cohomology circumstances we can give
a self-contained argument.}
 This is in
Section \ref{ev-coh}, where this machinery also provides a simple formula for
the Betti numbers of the polygon space (Section \ref {pp}).
 The form in which we
obtain the cohomology ring of the upper path space gives a simple guess for
the annihilator; to show this is the entire annihilator we use
a Gr\"obner basis argument, in Section \ref{72}.

The classical interest in these spaces is especially strong in the equal-weight
case, in that there is then an action of the symmetric group, and the quotient
is then a compactification of the moduli space of $m$-times-punctured Riemann
spheres. Our approach requires us to single out one edge, breaking this
symmetry.
However, a circle bundle associated to each edge gives a natural list
of degree 2 classes, permuted by the action of $S_m$ in
the equilateral case. In Section \ref{47} we locate these in our
presentation, and show they generate the ${\bbz}[{1\over 2}]$
cohomology ring. This gives a manifestly symmetric presentation
which is actually simpler.

But breaking the symmetry is unavoidable, in a very precise sense, if one
wants to compute the {\em integer} cohomology ring. While the action of
$S_m$ on the second rational cohomology group is the standard one on $\bbq^m$
\cite{Kl}, it is {\em not} the standard one on the second integral
cohomology group -- there is no $\bbz$-basis permuted by $S_m$.
In section \ref{equil} we show this, but also show that the action becomes
standard if one inverts (the necessarily odd number) $m$.

The main reason to avoid inverting 2 is to compute the ${\bbz}/2$-cohomology
ring of the {\em planar} polygon space, which we do in Section \ref{plpol}.

Lastly, if the edge chosen is the longest one, our formulae are no worse and
frequently much more computationally effective than the symmetric versions with
${\bbz}[{1\over 2}]$ coefficients. This and much else can be seen in Section
 \ref{Exples}
on examples.

{\em Acknowledgements.} Both authors are grateful to Shaun Martin for many
useful conversations concerning this problem. Susan Tolman made valuable
comments on a preliminary version of this paper. The second author wishes to
thank the University of Geneva for its hospitality while this work was being
done.

\vskip 1 truecm

\tableofcontents


\vskip 1 truecm

\section{The polygon spaces}  \label{definitions}

Let $\alpha =(\alpha _1,\dots ,\alpha _{m}) \in \rep^m$.
Let $S^2_{\alpha_i}$ denote the sphere in $\bbr^3$ with radius $\alpha_i$.
(Identifying $\bbr^3$ with $so(3)^*$, the Lie-Kirillov-Kostant-Souriau
symplectic structure gives $S^2_{\alpha_i}$ the symplectic volume $2\alpha_i$.)
Let us consider the manifold
$$ \prod_{i=1}^{m-1} S^2_{\alpha _i} \subset (\bbr^3)^m$$
equipped with the product symplectic structure.
We imagine $(\rho _1,\dots ,\rho _m)
\in (\bbr^3)^m$ as a path starting from the origin
of $m$ successive steps $\rho _i$ and thus call
 $\prod_{i=1}^{m} S^2_{\alpha _i}$  the {\it path space} for $\alpha $.
 The Hamiltonian actions on $\prod_{i=1}^{m} S^2_{\alpha _i}$
 which are relevant to us are
 \decale{a)} the diagonal $SO_3$-action with moment map
 $\mu(\rho):= \sum_{i=1}^{m-1} \rho _i$, ``endpoint"\smallskip
 \decale{b)} its restriction to $SO_2$ with moment map
 $\overline \mu (\rho ) :=  \zeta \big(\sum_{i=1}^{m-1} \rho _i\big)$,
  ``height of endpoint"
 \smallskip
 \decale{c)} the $(SO_2)^m$-action with moment map
 $\hat \mu (\rho ) = \big(\zeta (\rho _1),\dots ,\zeta (\rho _{m-1})\big)$,
 ``height of each step".

 This action makes
 $\prod_{i=1}^{m-1} S^2_{\alpha _i}$ a toric manifold. The moment
 polytope (image of $\hat \mu $) is the box $\prod_{i=1}^{m-1} [-\alpha _i,\alpha
 _i]$. 
 \vskip 0.5 truecm

  Let $\alpha =(\alpha _1,\dots ,\alpha _m) \in \rep^m$.
 We consider
 the space $\pol(\alpha)$ of configurations in $\bbr^3$
 of a polygon  with $m$ edges of length
 $\alpha _1,\dots ,\alpha _m$, modulo rotation, and call it
 the {\it polygon space} for $\alpha $.
 The precise definition is
 $$\pol(\alpha):=\bigg\{(\rho _1,\dots,\rho _m)\in (\bbr^3)^m \biggm |
 \forall i, |\rho _i|=\alpha_i
 \hbox{ and } \sum_{i=1}^m\rho _i=0  \bigg\}\biggm/SO_3$$
 where $SO_3$ acts on $(\bbr^3)^m$ diagonally.
 We say that $\alpha $ is {\it generic} if the equation
 $\sum_{i=1}^m\varepsilon _i\alpha  _i=0$ has no solution
 with $\varepsilon _i = \pm 1$ (there are no lined configurations). In this
 paper, $\alpha $ will always be assumed generic.

 The spaces $\pol (\alpha )$  have been
 studied for instance in \cite{Kl}, \cite{KM} and \cite{HK}.
 When $\alpha $ is generic it is shown that  $\pol (\alpha )$ is a
 closed smooth manifold of dimension $2(m-3)$ which naturally carries a
 symplectic form
 $\omega $. More precisely,
 $\pol (\alpha )$ occurs as symplectic reductions of the path spaces
 $$ \pol (\alpha ) = \mu^{-1} (0)\big/ SO_3 = \bigg(\prod_{i=1}^m S^2_{\alpha _i}\bigg)
 \mathop{\,\bigg/\!\!\!\bigg/\,}_{\!\textstyle 0\kern 7pt}SO_3
 \cong
 \bigg(\prod_{i=1}^{m-1} S^2_{\alpha _i}\bigg)
 \mathop{\,\bigg/\!\!\!\bigg/\,}_{\!\textstyle \alpha _m}SO_3 .$$
 This last object is ``paths of $m-1$ steps of lengths
 $\alpha_1,\ldots,\alpha_{m-1}$, whose endpoint is at distance $\alpha_m$
 from the origin, modulo $SO(3)$'', and obviously corresponds to the
 polygons as previously described. It is this last picture,
 in which the $m$th edge plays a distinguished role,
 that will be of most use to us.

 The {\it abelian polygon space} $\apol(\alpha )$ is defined as
 $$\apol(\alpha ):= \bigg\{(\rho _1,\dots,\rho _m)\in (\bbr^3)^m \biggm |
 |\rho _i|=\alpha _i \hbox{ and }
 \zeta\big(\sum_{i=1}^{m-1}\rho _i\big)=\alpha _m \bigg\}
 \biggm/SO_2.$$
 The word ``abelian" is used because $\apol(\alpha )$ is a symplectic
 reduction of the path space by the maximal torus $SO_2$ of $SO_3$:
 $$\apol(\alpha ) =\zeta^{-1}(\alpha _m)\big/SO_2= \prod_{i=1}^{m-1} S^2_{\alpha _i}
 \mathop{\,\bigg/\!\!\!\bigg/\,}_{\!\textstyle \alpha _m}SO_2.$$

 $\apol(\alpha )$ is visualized as the space of piecewise linear $(m-1)$-chains
 (with edge lengths $\alpha _1, \dots , \alpha _{m-1}$) which terminate on the
 plane $z=\alpha _m$, modulo rotations about the $z$-axis. The symplectic
 manifold $\apol(\alpha )$ is of dimension $2(m-2)$ and contains $\pol(\alpha )$
 (those that also end on the $z$-axis) as a symplectic submanifold of codimension
 $2$.

The $(SO_2)^{m-1}$-action on $\prod_{i=1}^{m-1} S^2_{\alpha _i}$ descends to a
Hamiltonian action on $\apol(\alpha )$. It is effective once we divide by the
diagonal subgroup $SO_2$ in $(SO_2)^{m-1}$. With this action,
$\apol(\alpha )$ is a toric manifold with moment polytope
$$\Xi_\alpha :=\bigg\{(x_1,\dots ,x_{m-1})\in
\prod_{i=1}^{m-1}[-\alpha _i,\alpha _i]
\biggm |  \sum_{i=1}^{m-1}x_i = \alpha _m\bigg\}.$$

 The {\it upper path space} $\upol(\alpha )$ is defined as
$$\upol(\alpha ):= \bigg\{\rho =(\rho _1,\dots,\rho _{m-1})\in
 (\bbr^3)^{m-1}\biggm |
\zeta(\sum_{i=1}^{m-1}\rho _i)\geq \alpha _m \hbox{ and } |\rho _i|=\alpha
_i\bigg\}\bigg/\sim$$
where the equivalence relation ``$\sim $" is defined as follows:
$\rho \sim \rho '$ if $\rho =\rho '$ or if
$$\zeta\big(\sum_{i=1}^{m-1}\rho _i\big)= \alpha _m \ \hbox{ and } [\rho ] =
 [\rho ']
\hbox{ in } \apol(\alpha ).$$
One can see $\upol (\alpha )$ as the result of a symplectic cut
(in the sense of
\cite{Le}) of the path space $\prod_{i=1}^{m-1} S^2_{\alpha _i}$
at the level $\alpha _m$ of the moment map $\overline \mu$. The space $\upol
(\alpha )$ is thus a closed symplectic  manifold of dimension $2(m-1)$ and is a
compactification of some open set of $\prod_{i=1}^{m-1} S^2_{\alpha _i}$. It
 contains $\apol (\alpha )$ as a codimension $2$ symplectic submanifold
 (one dimension lost by the height restriction, the other by the circle quotient).

The  Hamiltonian action of the torus  $(SO_2)^{m-1}$ on the path space descends
to an effective action on $\upol(\alpha )$. Therefore, $\upol (\alpha )$ is a
toric manifold with moment polytope $$\hat \Xi_\alpha := \bigm\{(x_1,\dots
,x_{m-1})\in \prod_{i=1}^{m-1}[-\alpha _i,\alpha _i] \ \bigm | \
\sum_{i=1}^{m-1}x_i \geq\alpha _m\bigm\}.$$ Observe that the codimension $2$
submanifold $\apol(\alpha ) \subset
 \upol(\alpha
)$ corresponds to the  facet ($=$ codimension 1 face) $\Xi_\alpha$ of the moment
polytope  $\hat\Xi_\alpha$ for $\upol (\alpha )$.  We shall prove now the
 important
fact that $\pol (\alpha )$ is obtainable as a transverse intersection of $\apol
(\alpha )$ with itself.

 The {\em vertical path space} is defined by
 $$\vpol (\alpha ):= \big \{\rho \in \upol (\alpha ) \bigm |
  \rho \hbox{ terminates on the $z$-axis }\big\}\subset \upol (\alpha )\, . $$

 It is a codimension 2 submanifold of  $\upol (\alpha)$ which we now show 
 intersects $\apol (\alpha )$ transversally in $\pol (\alpha )$.
 Consider the open subset  $W:=\mu^{-1}(\bbr^3-\{0\})$ of the path space
 $\prod_{i=1}^{m-1} S^2_{\alpha_i}$
 %
 The map  $\rho\mapsto \mu(\rho)/|\mu(\rho)|$
 is a fibration $W\to S^2$ (see \cite[(1.3)]{Ha}). The fiber
 $V$ over $(0,0,1)$ is an $SO_2$-invariant codimension 2 submanifold of $W$
 which project onto $\vpol (\alpha )$.
 The manifold $W$ intersects
 $M:=\hat \mu ^{-1}(\{\alpha _m\})$ transversally  in $\mu^{-1}(0,0,\alpha _m)$.
 The $SO_2$-action on $W$ induces a $SO_2$-action on $\upol (\alpha )$ for which
 $\apol (\alpha )$ is a set of fixed points.
 For each $\rho\in \pol (\alpha )$, the tangent space $T_\rho \upol (\alpha )$
 decomposes
 $$T_\rho \upol (\alpha ) \simeq T_\rho \apol (\alpha ) \oplus \bbr\, v \oplus
 \bbr\, r(v)$$
 where $v\in T_ \rho \vpol (\alpha )$ and $r\in SO_2$ is the rotation of angle
 $\pi/4$. As $\vpol (\alpha )$ is $SO_2$-invariant, the vector $r(v)$ belongs to
 $T_ \rho \vpol (\alpha )$ and thus
 $T_\rho \upol (\alpha ) \simeq T_\rho \apol (\alpha ) +  T_\rho \vpol (\alpha
 )$.

 \begin{Proposition} \label{isotopy}
 There is a smooth isotopy
 $\varphi _t: \vpol (\alpha ) \fl{}{} \upol (\alpha)$ such that
 $\varphi _0(\rho ) = \rho$ and
 $\varphi _1(\vpol (\alpha ))=\apol (\alpha )$.
 \end{Proposition}

 \preu $\varphi _t (\rho )$ is the image of $\rho $ by the rotation about the
 $y$-axis of angle
 $$t\,\cos^{-1}\frac{\alpha _m}{\zeta (\rho )}\ . $$

 \begin{Corollary} \label{tran-int}
 a) $\pol (\alpha )$ is a
 transverse intersection of $\apol (\alpha )$ with itself.
 \decale{b)} $\vpol (\alpha )$ and $\apol (\alpha )$ are diffeomorphic
 rel their common $\pol (\alpha)$.
 \end{Corollary}

 We note in passing that both $\upol(\alpha )$ and $\apol(\alpha )$
 are
 symplectomorphic to polygon spaces, though we will not use this fact
 elsewhere in the paper.

 \begin{Proposition} \label{longedges}
 For $\delta $ big enough, one has  symplectomorphisms
 $$\upol(\alpha _1,\dots ,\alpha _m) \cong \apol(\alpha _1,\dots ,\alpha _{m-1},
 \delta -\alpha _m ,\delta )$$
 $$\apol(\alpha _1,\dots ,\alpha _m) \cong \pol(\alpha _1,\dots ,\alpha _{m-1},
 \delta +\alpha _m, \delta ).$$
 \end{Proposition}

\preu  The symplectomorphisms come from the fact that the above spaces are toric
manifolds with isomorphic moment polytopes. For the first one, the moment
 polytopes
are
$\Xi_{\alpha _1,\dots ,\alpha _{m-1},\delta -\alpha _m ,\delta }$
 and $\hat\Xi_{(\alpha _1,\dots ,\alpha _m)}$ and
 $(x_1,\dots ,x_m)\mapsto
(x_1,\dots ,x_{m-1})$ gives the required isomorphism.
The condition on $\delta $ is $\delta >\alpha _m $.

For the second case, let
$\delta > \sum_{j=1}^{m-1}\alpha _j$. Consider the functions \newline
$d_i^{\ } : \pol(\alpha _1,\dots ,\alpha _{m-1},
\delta +\alpha _m, \delta ) \to \bbr$ ($i=0,\dots , m-1$) given by
$$d_i(\rho ) := \big\| -\rho _m +  \sum_{j=1}^{i}\alpha _j \big\|.$$
Because of the condition on $\delta $, the functions $d_i$ never vanish and are
 thus
smooth. It was shown in \cite{KM} that they Poisson-commute and generate an
 effective
$T^{m-2}$-action ($d_0=\delta $ and $d_{m-1}=\delta +\alpha _m$ are constant).
 This
makes  $\pol(\alpha _1,\dots ,\alpha _{m-1},
\delta +\alpha _m, \delta )$ a toric manifold with moment polytope $\Delta $
the set of $(x_0,\dots ,x_{m-1}) \in \bbr^{m}$ satisfying
$x_0=\delta$, $x_{m-1}=\delta  +\alpha _m $ and the triangle inequalities:
$$x_{i-1} + \alpha _i \geq x_i,\quad
x_i + \alpha _i \geq x_{i+1},\quad
x_{i-1} + x_i \geq  \alpha _i$$
for $i=1,\dots ,m-1$ (see \cite{HK}, \SS 5). By the condition on $\delta $, the
 third
inequality is automatically satisfied and the two others amount to $x_i \in
[x_{i-1}-\alpha _i ,x_{i-1}+\alpha _i]$ ($i=1,\dots ,m-2$). Changing
 variables
$$(x_0,\dots ,x_{m-1})  \mapsto \big(x_1-x_0,x_2-x_1,\dots
 ,x_{m-1}-x-{m-2}\big)$$
gives an isomorphism between $\Delta $ and $\Xi_{(\alpha _1,\dots ,\alpha
 _{m-1},
\delta +\alpha _m, \delta )}$.   \cqfd

\vskip .5 truecm

\section{Short and long subsets} \label{shortlong}

 Let $\alpha = (\alpha _1,\dots ,\alpha _m) \in \rep^m$. A
 subset $J \subset \{1,\dots ,m\}$ is {\it short} if
 $$ \sum_{j\in J} \alpha _j \leq \sum_{j\not\in J} \alpha _j\ .$$
  Equivalently, $J$ is short iff
 $$\sum_{j=1}^m(-1)^{{\textstyle \chi}
 _J(i)}\alpha _j \geq 0$$
 where $\chi _{\scriptscriptstyle S} : \bbn \to \{0,1\}$
 is the characteristic function of $S$. For example, the empty set is short, 
 and singletons are short
 iff $\pol (\alpha )\not = \emptyset$. More generally a set $S$ is short
 exactly if there exist configurations in $\pol (\alpha )$ with all edges in $S$
 parallel.
 Observe that the equalities in the above definitions cannot occur  if $\alpha $ is
 assumed to be generic. Define
 $$ \cals  := \cals (\alpha ):=  \big\{J \subset \{1,\dots ,m\} \bigm| J\hbox{ is short }\big\}
 \ .$$

 The collection $\cals$ is partially ordered by inclusion.
Every subset of a short subset is short,
and thus the poset $\cals$ is determined by its maximal elements.

  \begin{ex} \label{exquadri} \rm For $m = 4$, there are, up to poset isomorphism,
   two possibilities for $\cals$. Listing only the maximal subsets, they are:
  \decale{a)} $\cals (\alpha) \supset \big\{ \{1,2\},\{1,3\},\{2,3\},\{4\}\big\}$; example:
 $\alpha = (1,1,1,2)$.
   \decale{b)} $\cals (\alpha ) \supset \big\{ \{1,4\},\{2,4\},\{3,4\}\big\}$; example:
 $\alpha = (2,2,2,1)$.
 \end{ex}

 While the length vector $\alpha \in \rep^m$ defines  $\pol (\alpha )$ up to
 symplectomorphism, we shall see that the diffeomorphism type of $\pol (\alpha )$ is determined
 by the combinatorial data $\cals (\alpha )$.

 Define $F_\alpha : \prod_{i=1}^m S^{2} \fl{}{} \bbr^3$ by $F_\alpha
 (z_1,\dots ,z_m):= \sum_{i=1}^m\alpha _i z_i$.
 The map $F_\alpha $ is smooth and
 $SO_3$-equivariant for the diagonal action of $SO_3$ on
 $\prod_{i=1}^m S^{2}$ and the natural action on $\bbr^k$.
 Let $A_\alpha := F_\alpha^{-1}(0)$.
 By what is in Section \ref{definitions},
 the manifold $\pol (\alpha )$ is diffeomorphic to
 $A(\alpha )/SO_3$  and $A(\alpha )$ is the total space of a principal $SO_3$-bundle $\xi(\alpha)$.

  \begin{Proposition} \label{isodiffeo} Let $\alpha $ and $\alpha '$
  be generic elements in
 $\rep^m$. A poset isomorphism
 $\phi : \cals(\alpha )\fl{\simeq}{} \cals(\alpha ')$ determines
 (up to isotopy) an $SO_3$-equivariant diffeomorphism  between
 $ A (\alpha )$ and $A(\alpha ' )$.
 \end{Proposition}

 \preu The poset isomorphism $\phi$ is first of all a permutation of
 $\{1,\dots ,m\}$. The correspondence $\rho_i\mapsto \rho_{\phi (i)}$
 gives an equivariant
 diffeomorphism from $A (\alpha_1,\dots ,\alpha _m )$  onto
 $A (\alpha_{\phi (1)},\dots ,\alpha _{\phi (m)})$. Therefore, one can assume
 that $\cals (\alpha ) = \cals (\alpha ')$ and
 that $\phi = {\rm id}$.

 For $t\in[0,1]$, let $\alpha (t):= t\alpha + (1-t)\alpha ' \in \rep^m$.
 Define $\beta : [0,1]\times\prod_{i=1}^m S^{2}\to [0,1]\times\bbr^3$ by
 $\beta (t,z_1,\dots ,z_m):= (t,F_{\alpha (t)} (z_1,\dots ,z_m))$.
 The inequalities involved in the definition of $\cals(\alpha )$ are all strictly
 verified for $\alpha (t)$. This shows that $\alpha (t)$ is generic for all $t$.
 One deduces that all the points of
 $[0,1]\times \{0\}$ are regular values of $\beta$.
 The manifold $W:=\beta^{-1}([0,1]\times \{0\})$ is then an $SO_3$-equivariant cobordism from
 $M (\alpha )$ to $M(\alpha ')$ and $\beta{\mid W} : W \to [0,1]\times
 \{0\}$
 is an $SO_3$-invariant map without critical points. Choose an  $SO_3$-invariant
 Riemannian metric on $W$. The gradient flow of $\beta$ on $W$ then produces an $SO_3$-equivariant
 diffeomorphism from $A (\alpha )$ to $A({\alpha '})$.  \cqfd

 \begin{Corollary} \label{isodiffeocor}  If $\alpha $ and $\alpha '$ are generic and
 $\cals (\alpha )  \simeq \cals (\alpha ')$ then there is a diffeomorphism 
 $h : \pol (\alpha) \to \pol (\alpha ')$ such that
 $h^*\xi_{\alpha '}=\xi _\alpha $.
 \end{Corollary}

 \begin{Remarks}  \rm  a)\  For the  examples of \ref{exquadri}, both
  $\pol (1,1,1,2)$ and $\pol (2,2,2,1)$ are diffeomorphic to the sphere $S^2$.
  But $\xi(1,1,1,2)$ is the non-trivial $SO_3$-bundle over $S^2$ whereas
 $\xi(2,2,2,1)$ is the trivial one (see example \ref{quadr}).

  \decale{b)} We have no counterexample to the  converse of Corollary  \ref{isodiffeocor}.

  \decale{c)} Proposition \ref{isodiffeo} works for polygons in $\bbr^k$
 (see [Ha]). But, for $k>3$, even if
 $\alpha $ is generic, the action of $SO_k$ is not free on $A(\alpha )$
 and thus Corollary \ref{isodiffeocor} does not make sense.
\end{Remarks}

  Let $\alpha \in \rep^m$ and let $\cals:=\cals (\alpha )$. For $k\in
 \{1,2,\dots ,m\}$,
 we introduce the subposet $\cals_k $ of $\cals$:
 $$\cals_k = \cals_k(\alpha ):= \{J\subset \{1,\dots ,m\}-\{k\} \bigm | J\cup \{k\} \in \cals \}
 .$$
 In the subsequent sections, we give the Poincar\'e polynomial and
 presentations of the cohomology ring of $\pol (\alpha )$ in terms of
 $\cals_m$. Proposition \ref{link} below together with Corollary \ref{isodiffeocor} imply that the
 diffeomorphism type of $\pol (\alpha )$ is determined by
 any of the subposets $\cals_k$.

 \begin{Proposition} \label{link} Let $\alpha $ and $\alpha '$ be generic elements in
 $\rep^m$. Suppose that there are $k,k'\in \{1,2,\dots ,m\}$ such that there is
 a poset isomorphism
 $\varphi : \cals_k (\alpha )\fl{\simeq}{} \cals_k' (\alpha ')$. Then any
 bijection
 $\Phi : \{1,2,\dots ,m\}\fl{\simeq}{} \{1,2,\dots ,m\}$ which extends $\varphi$
 and satisfies $\Phi (k)=k'$ is a poset isomorphism from $\cals (\alpha )$
 onto
 $\cals (\alpha ')$.
 \end{Proposition}

 \preu Let $\Phi$ be a permutation of $\{1,2,\dots ,m\}$ as in the statement  \ref{link}. By renumbering the components of $\alpha $ using the permutation
 $\Phi$, one can assume that $k=k'$, $\cals_k=\cals_k'$ and $\Phi={\rm id}$. It
 then suffices to prove that $\cals_k=\cals'_{k'}$ implies $\cals=\cals'$.

 Let $J\subset \{1,2,\dots ,m\}$ and let $\bar J := \{1,2,\dots ,m\}-J$.
 If
 $k\in J$, then $J\in \cals$ iff $J-\{k\}\in \cals_k$. If $k\notin J$, then
 $k\in \bar J$ and $J\in \cals$ iff $\bar J - \{k\}\notin \cals_k$. This gives a
 procedure to decide whether or not $J\in \cals$ by only knowing $\cals_k$.
 Therefore $\cals_k$ determines $\cals$. \cqfd \smallskip

 A  set $J$ which is not short is called {\it long}. The following notation will
 be used
 $$ \call   :=   \big\{J \subset \{1,\dots ,m\} \bigm| J\hbox{ is long
  }\big\}$$
  and
  $$ \call_m  :=  \big\{J \subset \{1,\dots ,m-1\} \bigm| J\cup \{m\}
\hbox{ is long }\big\} \subset \call \ .$$

\goodbreak
\section{Pairs of even-cohomology manifolds} \label{ev-coh}

We will call a topological space $X$ an {\it even-cohomology space} if its
 cohomology
groups  $H^*(X;\bbz)$ vanish for $*$ odd. We write $H^*(X)$ for $H^*(X;\bbz)$.

\begin{Proposition} \label{ex-seq}
 Let $M$ be a closed oriented manifold of dimension
$n$ with $n$ even. Let $Q$ be a closed oriented submanifold of $M$ of
codimension $r$. Suppose that $Q$ and $M-Q$ are  even-cohomology spaces. Then
one has short exact sequences
\begin{equation}
0 \fl{}{} H_{n-*}(M-Q) \fl{}{} H^*(M)\fl{i^*}{}  H^*(Q)\fl{}{} 0
 \label{seq-one}
\end{equation}
and
\begin{equation}
0 \fl{}{}H^{*-r}(Q)\fl{}{} H^*(M)\fl{j^*}{}  H^*(M-Q)\fl{}{} 0 \label{seq-two}
\end{equation}
where $i^*$ and $j^*$ are the ring homomorphism induced by the inclusions. In
particular, $M$ is an even-cohomology space.
 \end{Proposition}

\preu Let $T$ be a closed tubular neighborhood of $Q$ in $M$.
Consider the long cohomology exact sequence of the pair $(M,T)$ $$\cdots \fl{}{}
H^{*-1}(T)\fl{}{} H^*(M,T)\fl{}{} H^*(M)\fl{}{} H^*(T)\fl{}{}
H^{*+1}(M,T)\fl{}{}\cdots.$$ One has $H^*(T) \cong H^*(Q)$, and excision and
Poincar\'e duality produce the isomorphisms $$H^*(M,T) \cong  H^*(M - {\rm
int\,} T, \partial T) \cong
 H_{n-*}(M - {\rm int\,} T) \cong H_{n-*}(M - Q)$$
which give sequence (\ref{seq-one}).

Sequence (\ref{seq-two}) comes from the cohomology exact sequence of the pair
$(M,M - {\rm int\,}T)$. Indeed, one has $H^*(M - {\rm int\,}T)\cong H^*(M - Q)$
and  the isomorphisms $$H^*(M,M - {\rm int\,} T) \cong H^*(T,\partial T)
\cong H^{*-r}(Q)$$ are given  by excision and the Thom isomorphism.
\cqfd

\begin{Corollary} \label{pp-mani}
 Let $M$  and $Q$ be as in \ref{ex-seq}. The Poincar\'e polynomials of
$\,M$, $Q$ are  calculable from the Poincar\'e polynomial of $\,M-Q$ by
$$(1-t^r)P_Q(t) = P_{M-Q}(t) - t^n P_{M-Q}(1/t)$$ and $$(1-t^r)P_M(t) =
P_{M-Q}(t) - t^{n+r} P_{M-Q}(1/t).$$
\end{Corollary}

\preu  The two exact sequences of \ref{ex-seq} give the equations
\[ \begin{array}{lcrcc}
 P_M(t)  & =  & P_Q(t) &  +  & t^n P_{M-Q}(1/t) \\
  P_M(t)  & =  & t^r P_Q(t) &  +  & P_{M-Q}(t)
\end{array}
\]
from which the equations of \ref{pp-mani} are deduced. \cqfd

By \ref{ex-seq}, $H^*(Q)$ is isomorphic to the quotient of $H^*(M)$ by the ideal
 $\ker i^*$. We shall use the following:

\begin{Proposition} \label{ker}
Let $M$  and $Q$ be as in \ref{ex-seq}. The kernel of $i^* :
H^*(M)\fl{}{} H^*(Q)$ is the annihilator of the cup product by the Poincar\'e
 dual
class of $Q$.
\end{Proposition}

\preu That the annihilator contains $\ker i^*$ is a very general fact,
as we now show. Let $[M]\in H_n(M)$ and $[Q]\in H_{n-r}(Q)$ be the
fundamental classes.
 The
 Poincar\'e dual $q\in H^r(M)$ of $Q$ is determined by the equation
 $q \cap [M] = i_*([Q])$.

 Let $a\in H^*(M)$. By standard properties of cup and cap products (see
 \cite{Sp}, Chapter 5 \SS 6), one has $$i_*\big(i^*(a) \cap q)\big) = a \cap
 i_*([Q]) = a\cap (q \cap [M]) = (a\cup q) \cap [M].$$ Therefore, if $i^*(a)=0$
  then $a\cup q = 0$.

  The reverse implication is true if $i_*$ is injective, which will follow from
   $Q$ and $M-Q$ being
 even-cohomology spaces. By the universal coefficient theorem,
 an even-cohomology space is an even-homology space and, as in the proof of
 \ref{ex-seq}, one
 gets the short exact sequence $$0 \fl{}{} H_{*}(Q) \fl{i_*}{} H_*(M)\fl{}{}
  H^{n-*}(M-Q)\fl{}{} 0.$$ Therefore $i_*$ is injective and $a\cup q=0$ implies
 that $i^*(a)=0$. \cqfd

\section{Poincar\'e polynomials of polygon spaces} \label{pp}

 Before working out the cohomology rings of the polygon spaces in the next
 section, we give here their Betti numbers, in the form of the Poincar\'e
 polynomial. These are easy to obtain from Corollary \ref{pp-mani}.
 Different formulae for the Poincar\'e polynomial of $\pol (\alpha )$ were
 already obtained in \cite{Kl}, \SS 2.2, by different methods, as well as the
 following lemma (\cite{Kl}, Corollary 2.2.2):

\begin{Lemma} \label{polevcoh}
 $\pol (\alpha )$ is an even-cohomology space.
\end{Lemma}

\preu Let us consider the diagonal-length function $\delta: \pol (\alpha
 )\fl{}{}
\bbr $ given by   $\delta (\rho ):=|\rho  _m - \rho  _{m-1}|$. It is smooth if
$\alpha _{m-1}\not = \alpha _m$ which can be assumed since
by Proposition \ref{isodiffeo}
changing the $\alpha _i$'s slightly
for a generic $\alpha $ does not change the diffeomorphism class of
 $\pol
(\alpha )$

 By \cite[Theorem 3.2]{Ha}, $\delta $ is a Morse-Bott function. The critical
 points are of even index and are isolated except possibly for the  two extrema.
 The pre-image $M_{\rm max}$ of the maximum is either a point or $\pol(\alpha
 _1,\dots ,\alpha _{m-2},\alpha _m +\alpha _{m-1})$. For the pre-image
 $M_{\rm min}$ of the
 minimum, there are three possibilities:
 \decale{--} one point
 \decale{--} $\pol(\alpha _1,\dots ,\alpha _{m-2},\alpha _m
 -\alpha _{m-1})$
 \decale{--} a 2-sphere bundle over
 $\pol(\alpha _1,\dots ,\alpha _{m-2},\alpha _m -\alpha _{m-1})$ (when the
 minimum is $0$).

 This enables us to prove Lemma \ref{polevcoh} by induction on $m$.
 A 2-sphere bundle over an even cohomology space is an even cohomology
 space using the Gysin sequence. Therefore, in all the cases
 $M_{\rm min}$ and therefore $\pol (\alpha ) - M_{\rm max}$ are
 even cohomology manifolds. If $M_{\rm max}$ is a point, we are done. Otherwise,
 we use Proposition \ref{ex-seq} to deduce that $\pol (\alpha )$ is an even
 cohomology manifold.
 \cqfd
 \smallskip

 We now use the inclusion $\pol (\alpha ) \subset \apol (\alpha )$ to obtain the
 Poincar\'e polynomials for the various polygon spaces. They are given in terms
 of the posets $\cals:=\cals(\alpha )$ or $\cals_m:=\cals_m(\alpha )$
 introduced in \SS \ref{shortlong}.

\begin{Proposition} \label{ev-coh2} \label{37}
The open manifolds $\apol (\alpha )-\pol (\alpha )$ and $\upol (\alpha )-\apol
(\alpha )$ are both even-cohomology spaces with the same Poincar\'e polynomial,
${\displaystyle\sum_{J\in \cals_m}^{\ }\!\!\! t^{2\,|J|}}$.
\plaint{\endproc \count37=\enno}
\end{Proposition}

\preu  Using \ref{tran-int}, one
 can replace $\apol (\alpha )$ by $\vpol (\alpha )$.
 Consider the function $d :\upol (\alpha ) - \apol (\alpha )$ defined by
$d(\rho):= -\zeta\big(\sum_{i=1}^{m-1}\rho _i\big)$ and denote by $d^v$ its
 restriction
to $\vpol (\alpha ) - \pol (\alpha )$.

By \cite{Ha}[Theorem 3.2] the map $d^v$
 is a Morse function with a critical
point of index $2|J|$ for each $J\in \cals_m$ (the critical point is the lined
configuration with all the $\rho _i$ pointing upwards if $i\notin S$ and
 downward
otherwise). This proves \ref{37} for $\apol (\alpha )$.

A small rotation around a horizontal axis will decrease
 $\zeta\big(\sum_{i=1}^{m-1}\rho_i\big)$ and so increase $d$. The slope is
 positive
away from $\vpol (\alpha )$ and thus $d$ has no  critical points other than
those of $d^v$. At one of these critical points, the rotation can be used to
 check
the  non-degeneracy and show that the index is the same for $d$ as for $d^v$.
 \cqfd

The  above two results,
using \ref{pp-mani}, give the following:

\begin{Corollary} \label{38}
 The various polygon spaces are even-cohomology spaces.
Their Poincar\'e polynomials are
$$P_{\pol (\alpha)}=
\frac{1}{1-t^2} \sum_{J\in\cals_m} (t^{2|J|}-t^{2(m-|J|-2)})$$
$$P_{\apol (\alpha)}=
\frac{1}{1-t^2} \sum_{J\in\cals_m} (t^{2|J|}-t^{2(m-|J|-1)})$$
$$P_{\upol (\alpha)}=
\frac{1}{1-t^2} \sum_{J\in\cals_m} (t^{2|J|}-t^{2(m-|J|)}).$$
\end{Corollary}

 \begin{Remark} \label{38b}  \rm The following expression for $P_{\pol (\alpha)}$
 was obtained, using another method, by Klyachko \cite[Theorem 2.2.4]{Kl}:
 $$P_{\pol (\alpha)}= \frac{1}{t^2(t^2-1)}\bigg(
 (1+t^2)^{m-1}-\sum_{J\in\cals}t^{2|J|}\bigg).$$
 This formula gives $P_{\pol (\alpha)}$ in terms of $\cals (\alpha )$ whereas
 those of \ref{38} are in terms of $\cals_m (\alpha )$. This illustrates that $\cals_m$
 determines $\cals$ (Proposition \ref{link}).
 Both expressions have advantages: the one in terms of
 $\cals$ is more symmetric whereas those using $\cals_m$
 have many fewer monomials and lead to easier  computations
(see \SS \ref{Exples}).
\end{Remark}

\begin{Corollary} \label{signatur}
If $m$ is odd, so that the dimension of $\pol(\alpha)$ is a multiple of 4,
the signature of the polygon space $\pol(\alpha)$ is
$\sum_{J\in\cals_m} (-1)^{|J|}$.
\end{Corollary}

\preu
In \cite{Kl} it is proven that the Hodge numbers $h^{pq}$ of the K\"ahler
manifold $\pol(\alpha)$ vanish except for the diagonal $h^{pp}$.
Then the Hodge signature theorem \cite{GrHa} implies that the
signature is the Poincar\'e polynomial evaluated at $t=i$.
\cqfd

This in turn is the Euler characteristic of the associated planar polygon
space (discussed further in section \cit{plpol}), and one plus that of the
poset $\cals_m - \{\emptyset\}$ which is a simplicial complex.

\vskip 0.3 truecm

\goodbreak
\section{The cohomology of the upper path space}
 \label{45}

In this section we give a presentation of the cohomology ring of the upper
path space,  the toric manifold with moment polytope $\hat\Xi_\alpha $. The
cohomology ring of a toric manifold  is given by Danilov's Theorem (see
\cite[Chapter 5]{Fu}, \cite[Theorem 4.14]{DJ}) which we recall below in a version
 useful for us.

Let $M^{2n}$ be a compact symplectic toric manifold
(acted on by the standard torus
$T^n= \bbr^n/\bbz^n$). Suppose that the moment polytope $\Delta :=\mu (M)
 \subset
\bbr^n$ is  given by a  family of inequalities indexed by a finite set $\calj$:
$$\Delta\, =\big\{x\in \bbr^n
\mid \langle x,w_j\rangle  \leq \lambda _j,\ j\in \calj \big\}$$
where $w_j\in \bbz^n$ is primitive and $\lambda _j\in\bbr $. Let $\calf _j$ be
 the
hyperplane $\{x\in \bbr^n
\mid \langle x,w_j\rangle  = \lambda _j\}$. We suppose that the $\calf_j$'s are
distinct. As $\calj$ is finite, the facet-hyperplanes of $\Delta $ must belong
 to the
family  and will be indexed by $\calj_0 \subset \calj$. Observe that
$$ j\in\calj_0 \quad \Leftrightarrow \quad
{\rm codim\,}\big(\calf_j \cap \Delta \big) = 1  .$$
Let $\{e_1,\dots ,e_n\}$ be a basis of $\bbr^n$. Danilov's Theorem gives a
presentation of the ring $H^*(M)$ with a generator $F_j\in H^2(M)$ for each
hyperplane $\calf_j$:

\plaint{\cleartabs \settabs
\+ \kern 1 truecm  &  $\sum_{j \in \calj} \langle e_i,v_j \rangle F_j$
\kern .7 truecm   &    if  \quad
${\rm codim\,} \bigcap_{j\in \calb}\,(\calf_j\cap \Delta ) > |\calb |$
 \kern .7truecm  & \cr}

\begin{Theorem}  \label{danilov}  \label{59}
 {\bf (Danilov)}
The cohomology ring
$H^*(M)$ is the quotient of the polynomial ring $\bbz[F_j;\, j\in \calj]$, where
 each $F_j$ is of
degree 2,  by the ideal $\cali$ generated by  the two families of
 relators:\smallskip

\begin{tabbing}
 \kern 1 truecm  \=  $\sum_{j \in \calj} \langle e_i,v_j \rangle F_j$
\kern .7 truecm  \=    if  \quad
${\rm codim\,} \bigcap_{j\in \calb}\,(\calf_j\cap \Delta ) > |\calb |$
 \kern .7truecm \= \kill

 \>  $\sum_{j \in \calj} \langle e_i,v_j \rangle F_j$
 \>    $i=1,\dots ,n$ \kern .7truecm  \> (linear relators) \\ \\
\>  $\prod_{j\in\calb }F_j$  \>    if  \quad
${\rm codim\,} \bigcap_{j\in \calb}\,(\calf_j\cap \Delta ) > |\calb |$
\>  (intersection monomials).

 \end{tabbing}
\end{Theorem}
\plaint{\endproc \count59=\enno}
\medskip

\begin{Remarks} \rm  \label{DRk}
  \begin{enumerate}

\item  The statements of  Danilov's Theorem in the literature are
only for
 $\calj = \calj_0$, but any generator $F_j$ for $j\notin \calj_0$ is in
$\cali $ using b) with $\calb = \{j\}$.

\item  \label{DR2}  When $j\in\calj_0$, the preimage $\mu ^{-1}(\calf_j )$ is a
 codimension 2
submanifold representing the Poincar\'e dual class of $F_j$. \medskip

\item \label{DR3}
The class $[\omega ] \in H^2(M;\bbr)$ of the symplectic form satisfies
$$[\omega ]\ =\ \sum_{j\in \calj}\lambda _j F_j$$
(See \cite[p. 132]{Gu}; the different sign comes
from the fact that our vectors $w_j$ are
pointing out of $\Delta $, contrarily to those in \cite{Gu}).

  \end{enumerate}
\end{Remarks}

We now apply \ref{59} to $\Delta = \hat\Xi_\alpha $. Let $e_1,\dots , e_{m-1}$
 be
the standard basis of $\bbr^{m-1}$. The polytope $\hat\Xi_\alpha$ is the subset
 of
$\bbr^{m-1}$ subject to the inequalities \smallskip

\begin{tabbing}

\kern 1.5 truecm \= \+ \kill
 $\langle x,e_i \rangle \leq \alpha _i$  \= \kern 1.8 truecm  \=
$j = 1,\dots m-1$ \\  \smallskip
$\langle x,-e_i \rangle \leq \alpha _i$ \>  \> $j= 1,\dots m-1$ \\ \smallskip
 $\langle x,-\sum_{i=1}^{m-1}e_i \rangle
\leq -\alpha _m$

\end{tabbing}

The relevant hyperplanes will be called
$$\calu_i :=  \{\langle x,e_i \rangle = \alpha _i\} \ ,\quad
 \calv_i :=  \{\langle x,e_i \rangle = -\alpha _i\} \ \hbox{ and }
\calr := \big\{\langle x,\sum_{i=1}^{m-1}e_i \rangle = \alpha _m  \big\}$$
with corresponding classes $U_i, V_i, R \in H^2(\upol (\alpha ))$.
Set
$$\tilde \calu_i := \hat \mu ^{-1} (\calu_i) \ , \
\tilde \calv_i := \hat \mu ^{-1} (\calv_i) \ , \
\tilde \calr := \hat \mu ^{-1} (\calr). $$
The first two are those polygons whose $i$th step points straight up, or
straight down; the third is the abelian polygon space.

If $A\subset \{1,\dots ,m-1\}$, define
$\tilde \calu_A := \bigcap_{i\in A}\tilde \calu_i^{\ }$ or
$\tilde \calv_A := \bigcap_{i\in A}\tilde \calv_i$.

\begin{Lemma} \label{61}
 Let $A,B\subset \{1,\dots ,m-1\}$ such that $A\cap B = \emptyset $. Then the
 image under $ \hat \mu : \upol (\alpha ) \fl{}{} \bbr$ of
 $\,\tilde \calu_A \cap \tilde \calv_B$
 is the interval
 $$ \hat \mu \big(\tilde \calu_A \cap \tilde \calv_B\big)\ =\
 \bigg[ - \sum_{i=1}^{m-1} (-1)^{\chi_A(i)}\,\alpha _i\,,\,
 \sum_{i=1}^{m-1}  (-1)^{\chi_B(i)}\,\alpha _i\,\bigg]\cap
 [\alpha _m,\infty).$$
 \plaint{\endproc \count61=\enno}
 \end{Lemma}

 \preu Recall that $\hat \mu (\rho)$ is the height of the endpoint of
 $\rho$. The highest it can get is when all edges point straight up (except
 those in $B$, required to point down); the lowest is when all edges
 point straight down (except for those in $A$) or at $z = \alpha_m$. \cqfd

We now work out
the presentation of $H^2(\upol (\alpha ))$ given by Danilov's Theorem with all
 the
generators $U_i, V_i$ and $R$.
Recall that $\call$ is the collection of long subsets of $\{1,\ldots,m\}$ and
$\call_m$ the collection of subsets $L \subseteq \{1,\ldots,m-1\}$ such that
$L \cup \{m\} \in \call$.

\begin{Proposition}  \label{40}  \label{}
 The ring $H^*(\upol (\alpha ))$ is the quotient of the
polynomial ring generated in degree 2 by the classes $R$, $U_i$ and  $V_i$
 ($i=1,\dots
,m-1$),  divided by the ideal generated by the following relators \smallskip
\begin{tabbing}

 \kern 2 truecm \=\+ \kill
(a)\quad  \= $U_i-V_i-R$    \= \kern 1 truecm \= $i = 1,\dots
,m-1$  \\
 \smallskip  (b)  \> $U_iV_i$        \>\> $i = 1,\dots ,m-1$  \\

(c) \> $\prod_{i\in L} V_i$  \>\> $L\subset \{1,\dots ,m-1\}$ and $L\in \call_m$
\\ \smallskip
 (d)  \> $R\,\prod_{i\in L} U_i$  \>\> $L\subset \{1,\dots ,m-1\}$ and
$L\in \call$
\end{tabbing}
\plaint{\endproc \count40=\enno}
\end{Proposition}

\preu  The relators (a) are the linear relators of Danilov's
theorem. Clearly $\calu_i\cap \calv_i = \emptyset$ (an edge cannot point both up
and down) whence relators (b).

Suppose that $\rho \in \calv_L$. If
 $L\in \call_m$ then $\sum_{i=1}^{m-1}(-1)^{\chi_L(i)}\,\alpha _i < \alpha _m$.
 By
Lemma \cit{61},
 one has
$$\hat\mu (\rho ) \ \leq \
\sum_{i=1}^{m-1} (-1)^{\chi_B(i)}\,\alpha _i < \alpha _m$$
which contradicts $\rho \in \upol (\alpha )$. Therefore $\calv_L = \emptyset$ if
$L\in \call_m$ which gives relators (c). In words, a path that steps down too
much cannot end above $z\geq \alpha _m$.

Similarly, if $L \subset \{1,\dots ,m-1\}$  and
$\rho \in \calu_A$, then
$$\alpha _m < - \sum_{i=1}^{m-1} (-1)^{\chi_A(i)}\,\alpha _i \ \leq \  \hat\mu
 (\rho )$$
 and thus $\calr \cap \calu_A = \emptyset$ which produces relators (d).
 In words, a path that steps up
too much cannot end at $z= \alpha _m$.

We have thus proved that the families (b)--(d)
are indeed intersection monomials. We now prove that any intersection monomial
is a multiple of these. By Danilov's theorem, an intersection monomial
$C$ is square-free so of the form
$$C= \prod_{i\in A}U_i\,\prod_{j\in B}V_j \quad \hbox{ or } \quad
C= R\, \prod_{i\in A}U_i\,\prod_{j\in B}V_j.$$
If $i\in A\cap B $ then $C$ is a multiple of $U_iV_i$. Therefore, one may
 suppose that
$A\cap B = \emptyset$.

If $C=\prod_{i\in A}U_i\,\prod_{j\in B}V_j$ is an intersection monomial, then
${\rm codim\,} (\calu_A \cap \calv_B \cap \hat\Xi_\alpha ) \geq |A\cup B|$.
So by Lemma \cit{61} we know
$\sum_{i=1}^{m-1} (-1)^{\chi_B(i)}\,\alpha _i \leq \alpha _m$.
But genericity implies that this inequality is strict.
Therefore, the inequality 
 $\sum_{i=1}^{m-1} (-1)^{\chi_B(i)}\,\alpha _i < \alpha _m$ holds, that is
$B \in \call_m$ and thus  $C$ is a multiple of $\prod_{j\in B}V_j$,
an intersection monomial in (c).

Consider now the case
$C=R\, \prod_{i\in A}U_i\,\prod_{j\in B}V_j$. Thus
${\rm codim\,} (\calu_A \cap \calv_B \cap \calr \cap \hat\Xi_\alpha ) \geq
 |A\cup
B|+1$. By Lemma \cit{61}, this would not happen if
$$- \sum_{i=1}^{m-1} (-1)^{\chi_A(i)}\,\alpha _i \,<\,  \alpha _m  \,<\,
\sum_{i=1}^{m-1} (-1)^{\chi_B(i)}\,\alpha _i$$
and equalities never occur. We saw before that the inequality \\
$\sum_{i=1}^{m-1} (-1)^{\chi_B(i)}\,\alpha _i < \alpha _m$ makes
$C$ a multiple of an intersection monomial of (c). The other possibility is
$\alpha _m < - \sum_{i=1}^{m-1} (-1)^{\chi_A(i)}\,\alpha _i$ which is equivalent
 to
$A\in \call$ and makes $C$ a multiple of the intersection monomial $R\prod_{i\in
A}U_i$ of (d).  \cqfd  \smallskip

The presentation of $H^*(\upol (\alpha ))$  which will turn out to be useful is
the following one:

\begin{Theorem}  \label{41}
 The ring $H^*(\upol (\alpha ))$ is the quotient of the
 polynomial ring generated in degree 2 by the classes $R$ and  $V_i$ for
 $1\leq i\leq m-1$, divided by the ideal $\cali$ generated by the following
 families
 of relators \smallskip

 \begin{tabbing} \renewcommand{\arraystretch}{5}
\kern .7 truecm \= \+ (R1)\quad   \=
$R^2\sum_{S\subset L \atop S\in \cals_m} \big(\prod_{i\in S}V_i\big)
R^{|L-S|-1}$  \=    \= \kill

(R1)\> $V_i^2+ RV_i$ \>  \> $1\leq i\leq m-1$ \\ \bigskip

 (R2)  \> ${\displaystyle \prod_{i\in L}^{\ } V_i}$  \>\>
  $L\in \call_m$ \\ \bigskip

 (R3)  \> ${\displaystyle R^2\sum_{S\subset L \atop S\in \cals_m}^{\ }
 \big(\prod_{i\in
S}V_i\big) R^{|L-S|-1}}$  \>\> $L\subset \{1,\dots ,m-1\}$ and $L\in \call$
\end{tabbing}
\end{Theorem}

\preu  This presentation is obtained by algebraically transforming that of
 Proposition \cit{40}. The linear relations $U_i= V_i+R$ of \cit{40} are absorbed
 by
 reducing the family of  generators to $R$ and the $V_i$'s. Replacing $U_i$ by
 $V_i+R$ in monomials (b) gives  relators (R1). Relators  (R2)
 are just relators (c) (one could of course restrict to minimal
 sets $L\in \call_m$).
 Relators (R3) are obtained by expanding monomials (d):
 $$R\prod_{i\in L} U_i\  =  \ {\displaystyle R\sum_{S\subset L}
 \big(\prod_{i\in S}V_i\big) R^{|L-S|}}
  \  = \ {\displaystyle R\sum_{S\subset L \atop S\in
 \cals_m} \big(\prod_{i\in S}V_i\big) R^{|L-S|}}$$ (the second equality  is
 obtained thanks to relators (R2) which kill $\Pi_S V_i$ if $S\notin \cals_m$).
 As $L\in \call$ and
 $S\in \cals_m$, one has $S\not = L$. Therefore $|L-S| \geq 1$ and one can pull out
 one more $R$ to get relators (R3).   \cqfd

\vskip 0.5 truecm

\section{The cohomology rings of $\apol (\alpha )$ and $\pol (\alpha )$}
 \label{72}

\begin{Theorem}  \label{64}
The cohomology rings of the abelian polygon space and actual
polygon space can be obtained from that of the upper path space as follows:
$$ H^*( APol(\alpha)) \iso H^*( UPol(\alpha)) / \Ann(R) $$
$$ H^*(  Pol(\alpha)) \iso H^*( UPol(\alpha)) / \Ann(R^2) $$
where $\Ann(x) := \{ y \in H^*(UPol(\alpha)) : yx = 0 \}$.
\end{Theorem}

 \preu  By construction, $R$ is Poincar\'e dual in $\upol (\alpha )$ to $\apol (\alpha )$. By
 Proposition \ref{tran-int}, the Poincar\'e dual to $\pol (\alpha )$ in
 $\upol (\alpha )$ is $R^2$. All the spaces under consideration are
 even-cohomology spaces by \ref{polevcoh} and \ref{37}.
 Therefore the result follows from Proposition \ref{ker}. \cqfd

It remains to calculate these annihilators, or equivalently the
``ideal quotients''
$${\cal I}:R^k := \{ y \in \bbz[V_i,R] : R^k y \in {\cal I} \}, \qquad k=1,2$$
where ${\cal I}$ is the ideal found in Theorem \cit{41} defining
$H^*(UPol(\alpha))$. Manifestly these contain the families (1), (2),
and $R^{-k}$(3) (recall that $R^2$ divides all the relators in the
third family). We will show that these do in fact generate
the ideal quotients.

If (1)-(3) were a {\it Gr\"obner basis} for the ideal, this would
be straightforward (see the lemma below); it is not in general,
but we will show that it is close enough.

We take the computational viewpoint of Gr\"obner bases, that they provide
a recognition algorithm for elements of an ideal -- a polynomial is an
element of $\cal I$ if the reduction algorithm (defined below) can
reduce it to zero. Conveniently, any list of generators of an ideal
can be finitely extended to a Gr\"obner basis by adding S-polynomials
(also defined below);
if all S-polynomials reduce to 0, the basis is Gr\"obner.
While all necessary definitions are given here,
our reference for these theorems is \cite{Ei}.

Given a polynomial $p$ we wish to check for $\cal I$-membership,
a well-ordering of all monomials respecting multiplication ($a<b$ implies
$ac<bc$ for all $a,b,c$), and a list $\{r_i\}$ of generators of the ideal,
the reduction algorithm is defined as follows. Within each $r_i$ is an
initial monomial $m_i$ (with respect to the well-ordering). If one of those
$m_i$ divides a monomial $m_i l$ of $p$, ``reduce'' $p$ to
$p - r_i l$ (which is in $\cal I$ exactly if $p$ itself was).
This kills the $m_i l$ in $p$.
This algorithm terminates; $\{r_i\}$ is called a Gr\"obner
basis if $p\in {\cal I}$ implies that it terminates at zero,
{\em no matter what order the reductions take place}.
This powerful independence makes it very easy to prove things
about Gr\"obner bases.

Given two relations $r_1,r_2$ with initial monomials $m_1,m_2$, there
are two ways to reduce the monomial $lcm(m_1,m_2)$:
to ${m_2 \over gcd(m_1,m_2)} (m_1-r_1)$ and to
${m_1 \over gcd(m_1,m_2)} (m_2-r_2)$.
Their difference
$$ S(r_1,r_2) := {m_2 \over gcd(m_1,m_2)} (r_1-m_1)
                -  {m_1 \over gcd(m_1,m_2)} (r_2-m_2)$$
is called the {\it S-polynomial} of $r_1$ and $r_2$, and is
manifestly in the ideal; if the list $\{r_i\}$ cannot reduce these to 0,
it certainly isn't Gr\"obner.
There are two convenient converses to this fact \cite{Ei}:
\decale{(1)} if all the S-polynomials do reduce to 0, the list $\{r_i\}$ is a
 Gr\"obner
basis;
\decale{(2)} if not, one can add those S-polynomials as new elements of the
 list,
a process that eventually terminates at a Gr\"obner basis.

The following lemma points out the relevance of Gr\"obner bases to
calculating ${\cal I}:R^k$. In what follows we use an
{\it reverse lexicographic order} for $R$;
this means that monomials are ordered
first by their power of $R$ (with low powers earlier in the order),
and only then by other criteria (which we leave unspecified).

\begin{Lemma}  \label{L1}
Let ${\cal I} \leq \bbz[{x_i},R]$ be a homogeneous ideal in a polynomial ring,
and $\{r_i\}$ be a Gr\"obner basis of $\cal I$,
with respect to a ``revlex'' order for $R$.
Then $\{ r_i/gcd(r_i,R^k) \}$ is a Gr\"obner basis for ${\cal I}:R^k$.
\end{Lemma}
\plaint{\endproc}

\preu
First we show that each $r_i/gcd(r_i,R^k)$ is in fact in ${\cal I}:R^k$:
$$ R^k r_i\big/gcd(r_i,R^k) = {R^k \over gcd(r_i,R^k)} r_i \in {\cal I}.$$
Second, that this list $\{ r_i/gcd(r_i,R^k) \}$ is powerful enough to
reduce any element $p$ of ${\cal I}:R^k$ to zero. To see this, follow the
reductions of $R^k p$, an element of $\cal I$,
by the (assumed) Gr\"obner basis $\{r_i\}$.
The possible reductions of $R^k p$ using $\{r_i\}$ correspond  exactly
to possible reductions of $p$ using $\{r_i/gcd(r_i,R^k)\}$, because
reducing $R^k p$ using $r_i$ necessarily uses a multiple of $lcm(R^k,r_i)$,
which we divide by $R^k$ to get the corresponding reduction of $p$
using $lcm(R^k,r_i)/R^k = r_i/gcd(r_i,R^k)$.  \cqfd

This technique of linking one reduction algorithm to another will be
used again in what follows; we will say that the reductions are
{\it parallel} in $(p,\{r_i\})$ and
$(p',\{r'_i\})$ given a correspondence between possible
reductions of $p$ using $\{x_i\}$
and possible reductions of $p'$ using $\{x'_i\}$.

Unfortunately, the list of relations $(R1)-(R3)$ in Theorem \ref{41}
is not generally
a Gr\"obner basis, and extending it to one seems difficult -- in particular,
defining the problem would require more precise specification of the
monomial order, such as an ordering on the edges.
Luckily, this list is close enough to being Gr\"obner to calculate
the annihilators we need.

\begin{Theorem}  \label{groeb}
Let ${\cal I} \leq \bbz[\{x_i\},R]$ be a homogeneous ideal
in a polynomial ring, and $\{r_i\}$ generate $\cal I$ as an ideal.
Assume that all S-polynomials of pairs $r_i,r_j$, such that neither is
a multiple of $R^k$, reduce to zero
with respect to an elimination order for $R$.
Then $\{ r_i/gcd(r_i,R^k) \}$ generates ${\cal I}:R^k$ as an ideal.
\end{Theorem}

This is a weaker requirement than in Lemma \ref{L1},
which required that {\em all} S-polynomials reduce to zero;
this will let us ignore the third family of relators in \cit{41}.

\preu
The argument is this: we complete $\{r_i\}$ to a Gr\"obner basis for $\cal I$,
and show that the parallel completion of
$\{ r_i/gcd(r_i,R^k) \}$ is to a Gr\"obner basis of ${\cal I} : R^k$.
Therefore $\{ r_i/gcd(r_i,R^k) \}$ generates ${\cal I} : R^k$.

Let $r_1, r_2$ be generators such that $r_2$ is divisible by $R^k$.
Consider the S-polynomial $s := S(r_1,r_2)$.
We claim that the reductions are parallel for
$(s, \{r_i\})$
and
$(s/R^k, \{r_i/gcd(r_i,R^k)\}).$

\renewcommand{\div}{{\big |}}
For this to make sense, we first must establish that $R^k \div s$.
Let $R^j$ be the highest power of $R$ dividing $r_1$.
Then by our assumption on the order, $R^j \div m_1$ and $R^k \div m_2$.
In our formula for
$$ s := {m_2 \over gcd(m_1,m_2)} (m_1-r_1)
      - {m_1 \over gcd(m_1,m_2)} (m_2-r_2)$$
we can then see that $R^{k-j} \div  {m_2 \over gcd(m_1,m_2)}$,
$R^j \div (m_1-r_1)$, and $R^k \div (m_2-r_2)$, so $R^k \div S$.

Second, we must establish a correspondence between the possible reductions.
This is as before: reducing $s$ by adding a multiple of $r_i$ necessarily
adds a multiple of $lcm(r_i,R^k)$, which corresponds to reducing $s/R^k$
by adding a multiple of $lcm(r_i,R^k)/R^k = r_i/gcd(r_i,R^k)$.

Now consider the process of extending $\{r_i\}$ to a Gr\"obner basis
by tossing in a S-polynomial which cannot reduce to zero.
By the assumption, it must be of the above type
(one of the relations is divisible by $R^k$),
at which point it parallels an S-polynomial in $\{r_i/gcd(r_i,R^k) \}$.
What this establishes is that the $\{r_i/gcd(r_i,R^k) \}$
can generate a Gr\"obner basis of ${\cal I}:R^k$.
In particular they generate ${\cal I}:R^k$.
\cqfd

In the case at hand, with the three families of relations

 \begin{tabbing} \renewcommand{\arraystretch}{5}
\kern .7 truecm \= \+ (R1)\quad   \=
$R^2\sum_{S\subset L \atop S\in \cals_m} \big(\prod_{i\in S}V_i\big)
R^{|L-S|-1}$  \=    \= \kill

(R1)\> $V_i^2+ RV_i$ \>  \> $1\leq i\leq m$ \\ \bigskip

(R2)  \> ${\displaystyle \prod^{\ }_{i\in L \atop \{i\}\in \cals_m} V_i}$  \>\>
  $L\in \call_m$ \\ \bigskip

 (R3)  \> ${\displaystyle R^2\sum^{\ }_{S\subset L \atop S\in \cals_m}
 \big(\prod_{i\in
S}V_i\big) R^{|L-S|-1}}$  \>\> $L\subset \{1,\dots ,m-1\}$ and $L\in \call$
\end{tabbing}
we have to check the (R1)-(R1), (R1)-(R2), and (R2)-(R2) S-polynomials -- the
lemma lets us ignore the S-polynomials with (R3), since those are all divisible
by $R^2$.

One standard observation \cite{Ei} is that the S-polynomial of
two elements $(r_1,r_2)$ with relatively prime initial terms $m_1$, $m_2$
is necessarily trivial. Here the initial terms are
(R1) $V_i^2$  for (R1) and
 $\prod_{j\in L} V_j$ for (R2).

(R1)--(R1):\ Each pair of initial terms is relatively prime.

(R1)--(R2):\ $V_i^2$ and $\prod_{j\in L} V_j$ have greatest common divisor $V_i$
 if
$i\in L$, and are otherwise relatively prime.
$$ S(V_i^2 + R V_i, \prod_{j\in L} V_j)=
\big(\prod_{j\in L-\{i\}} V_j\big) (R V_i) - 0
                                = R \prod_{j\in L} V_j \equiv 0. $$

(R2)--(R2):\ The S-polynomial of two monomial relations is automatically zero,
no reduction necessary.

Using \cit{41} and \cit{64}, we have just proved

\begin{Theorem} \label{62}

 1)\ The cohomology ring of the abelian polygon space $\apol (\alpha)$ is $$
 H^*(\apol(\alpha)) =   \bbz[R,V_1,\dots ,V_{m-1}] / {\cal I}_{\apol}$$
 where $R$ and $V_i$ are of degree 2 and
  ${\cal I}_{\apol }$ is the ideal generated by the three families\smallskip
  \begin{tabbing} \renewcommand{\arraystretch}{5}
 \kern .7 truecm \= \+ (R1)\quad   \=
 $R^2\sum_{S\subset L \atop S\in \cals_m} \big(\prod_{i\in S}V_i )
 R^{|L-S|-1}$  \=    \= \kill

(R1)\> $V_i^2+ RV_i$ \>  \> $i=1,\dots ,m-1$ \\ \bigskip

(R2)  \> ${\displaystyle \prod^{\ }_{i\in L } V_i}$  \>\>
 $L\in \call_m$ \\ \bigskip

 (R3)  \> ${\displaystyle R\sum^{\ }_{S\subset L \atop S\in \cals_m} \big(\prod_{i\in
S}V_i\big) R^{|L-S|-1}}$  \>\> $L\subset \{1,\dots ,m-1\}$ and $L\in \call$
\end{tabbing}

 2) The cohomology ring of the actual polygon space $\pol(\alpha)$ is $$
 H^*(\pol) = \bbz[R,V_1,\dots ,V_{m-1}] / {\cal I}_{\pol} $$
 where $R$ and $V_i$ are of degree 2
 and ${\cal I}_{\pol}$ is generated by the three families\smallskip

 \begin{tabbing} \renewcommand{\arraystretch}{5}
\kern .7 truecm \= \+ (R1)\quad   \=
$R^2\sum_{S\subset L \atop S\in \cals_m} \big(\prod_{i\in S}V_i\big)
R^{|L-S|-1}$  \=    \= \kill

(R1)\> $V_i^2+ RV_i$ \>  \> $i=1,\dots ,m-1$ \\ \bigskip

(R2)  \> ${\displaystyle \prod^{\ }_{i\in L} V_i}$  \>\>
 $L\in \call_m$ \\ \bigskip

 (R3)  \> ${\displaystyle \sum^{\ }_{S\subset L \atop S\in \cals_m} \big(\prod_{i\in
S}V_i\big) R^{|L-S|-1}}$  \>\> $L\subset \{1,\dots ,m-1\}$ and $L\in \call$
\end{tabbing}
\end{Theorem}

As an example, we give the  expression of the class  $[\omega ] \in H^2(\pol
(\alpha ));\bbr ) $ of the symplectic form in terms of the generators $R$ and
$V_i$:

\begin{Proposition} \label{52}
The class $[\omega ] \in H^2(\pol (\alpha ));\bbr ) $  is given by
$$ [\omega ]\  =\
\bigg(-\alpha _m + \sum_{j=1}^{m-1}\alpha _j\bigg)\, R\ +\  2\sum_{\{j\}\in
\cals_m}\!\alpha _j V_j .$$
\end{Proposition}

\preu  From Remark \ref{DR3} of \ref{DRk}, one gets
$$ [\omega ]\  =\ -\alpha_mR + \sum_{i=1}^{m-1}\alpha_i\big(U_i - V_i\big)$$
which is put in the required form by using the relations $U_i=V_i+R$ of
\ref{40}. Observe that the formula of \ref{52} is actually also valid in $
H^2(\upol (\alpha );\bbr ) $ and in $H^2(\apol (\alpha );\bbr ) $. \cqfd

 As a consequence, one has a sufficient condition for the class
 $[\omega ] \in H^2(\pol (\alpha ));\bbr ) $ to be integral:

\begin{Corollary}  If  $\alpha \in \bbz^m$ then
$[\omega ] \in H^2(\pol (\alpha); \bbz ) $.
\end{Corollary}

\vskip .5 truecm

\section{Natural bundles over polygon spaces}  \label{47}

Let $\alpha =(\alpha _1,\dots ,\alpha _m)\in \rep^m$ be generic.  For
$j\in\{1,\dots ,m\}$ we define $A_j:=A_j(\alpha ) \subset (\bbr^3)^m$ by
$$A_j:=\bigg\{\rho =(\rho _1,\dots,\rho _m)\in (\bbr^3)^m\bigg| |\rho _i|=\alpha
_i \hbox{ and }
\sum_{i=1}^m\rho _i=0 \hbox{ and }
 \rho_j=(0,0,\alpha _j)\bigg\}.$$

As $\alpha $ is generic, the diagonal $SO_2$-action on $(\bbr^3)^m$ is free on
 $A_j$
and  one has $\pol(\alpha ) = A_j/SO_2$. Therefore,
$A_j\to
\pol(\alpha )$ is a principal $SO_2$-bundle $\xi _j$ determined by its Chern
(or Euler) class
$c_j := c_1(\xi _j) \in H^2(\pol(\alpha );\bbz)$.

As in Section \ref{shortlong}, let us consider
$$A:=A(\alpha ):=\bigg\{(\rho _1,\dots,\rho _m)\in (\bbr^3)^m\bigg|
\sum_{i=1}^m\rho _i=0 \hbox{ and } |\rho _i|=\alpha _i \bigg\} \subset
 (\bbr^3)^m.$$
 As $\alpha $ is generic, the quotient map $A\fl{}{} \pol(\alpha )$ is a
 principal $SO_3$-bundle denoted by $\xi :=\xi (\alpha)$
 (write the elements of
 $\bbr^3$ as row vectors, so that $SO_3$ acts on the right on them).
 The bundle $\xi$ is determined by
 its Stiefel-Whitney class $w_2(\xi )\in H^2(\pol(\alpha );\bbz_2)$
 and its Pontrjagin class $p := p_1(\xi ) \in H^4(\pol(\alpha );\bbz)$.

\begin{Proposition} \label{48}
For each $j\in\{1,\dots ,m\}$ the bundle $\xi _j$ is a
$SO_2$-reduction of $\xi $, that is, $A$ is $SO_3$-equivariantly
diffeomorphic to $A_j\, \times _{SO_2}\, SO_3$.
\end{Proposition}

\preu  The $SO_3$-equivariant diffeomorphism from
$A_j\, \times _{SO_2}\, SO_3$ onto
$A$ is given by $(\rho , \beta ) \mapsto(\rho )\beta $. \cqfd

\begin{Corollary} \label{49}
For each $j\in\{1,\dots ,m\}$ one has
$c_j^2=p$ and \hfb $c_j=w_2(\xi ) \ {\rm mod}\, 2$.
\end{Corollary}

\preu  The first equation is shown in \cite{MS}, Corollary 15.8, p 179.
The second is classical between Euler and Stiefel-Whitney
classes (\cite{MS}, Property 9.5, p. 99).
\cqfd

 By  \ref{62} the Chern classes $c_j$ are
expressible in terms of the classes
$R$ and $V_j$. The formulae are:

\begin{Proposition} \label{50}
In $H^2(\pol (\alpha);\bbz )$, one has \smallskip
$$ c_i = \left\{
\begin{array}{ll}
 R +2 V_i & \hbox{for   $i = 1,\ldots,m-1$} \\
  -R &  \hbox{for $i=m$}
\end{array} \right. $$
(In particular $c_i = R$ if $\{i\} \in \call_m$.)
\end{Proposition}

\preu
Define $B \subset (\bbr^3)^{m-1}$ by
$$B:=\big\{(\rho _1,\dots,\rho _{m-1})\in (\bbr^3)^m \bigm |
|\rho _i|=\alpha _i \hbox{ and }
\zeta\big(\sum_{i=1}^m\rho _i)=\alpha _m\big\}.$$
As $\alpha $ is generic, $SO_2$ acts freely on $B$ making it a principal
$SO_2$-bundle $\psi$ over  $\apol (\alpha )$. One has a commutative diagram
$$
\begin{array}{ccc}
A_m & \fl{\tilde i}{}  & B \cr
\downarrow   &  & \downarrow \cr
\pol (\alpha ) & \fl{i}{}  & \apol (\alpha )
\end{array}
$$
where the inclusion $\tilde i: A_m \hookrightarrow B$ is anti-equivariant:
 $j(\rho
\cdot \beta ) = \rho \cdot \beta ^{-1}$ (since, for the identification of $\pol
(\alpha )$ as a subspace of $\apol (\alpha )$, the vector $\rho _m$ must point
downwards).   Therefore $c_m = -i^*(c_1(\psi ))$ and the equality
$c_m= -R$ in $H^2(\pol (\alpha )))$ is equivalent to
$c_1(\psi)= R$ in $H^2(\apol (\alpha ))$ (recall that $R$ denotes a class in
$H^2(\upol (\alpha ))$ as well as its  its images in $H^2(\apol (\alpha ))$ and
$H^2(\pol (\alpha ))$).

Let $T$ be a tubular neighbourhood of $\tilde \calr = \apol (\alpha )$ in $\upol
(\alpha )$.  The retraction $T\to \tilde \calr $
is the disc bundle associated to $\psi$. The class
$R\in  H^*(\upol (\alpha ))$ being the Poincar\'e dual of $\tilde \calr $,
it is the
image of the Thom class,  ${\rm Thom}(\psi )\in H^2(T,\partial T)$, under the
homomorphism
$$H^2(T,\partial T) \fl{\cong}{} H^2(\upol (\alpha ), \upol (\alpha ) -
{\rm int\,}T) \fl{}{} H^2(\upol (\alpha )) .$$
Therefore $R\in H^2(\apol (\alpha ))$ is the image
of ${\rm Thom}(\psi )$ under the homomorphism 
$H^2(T,\partial T) \to H^2(T) \cong H^2(\upol (\alpha ))$ which, by
one of the definitions of the Euler class  (\cite{Hu}, \SS
16.7), is equal to $c_1(\psi)$. Thus, we have proven that $c_m=-R$.

 By the Duistermaat-Heckmann theorem \cite[Theorem 2.7]{Gu}, one has
 $$-c_1(\psi ) = \frac{\partial }{\partial \alpha _m}[\omega]$$
 \setcounter{equation}{0}
 in $H^2(\apol(\alpha );\bbr)$
 and thus
 \begin{equation}\label{eqcm} c_m = \frac{\partial }{\partial \alpha _m}[\omega]
 \end{equation}  in
 $H^2(\pol(\alpha );\bbr)$. Finally, by symmetry (any edge can be the ``last''
 one): $$c_j = \frac{\partial }{\partial \alpha _j}[\omega] \label{der}.$$

 Applying this formula to the expression of $[\omega ]$ given in Proposition
 \ref{52} gives the equations of \ref{50}.

\begin{Corollary} \label{51}
 The classes $c_i$ generate $H^*(\pol (\alpha );\bbz [1/2] )$.
\end{Corollary}

\begin{Remarks} \label{cherncl} \rm \ a) By \ref{49} the classes $c_i$ generate a 1-dimensional
space in \\ $H^2(\pol(\alpha );\bbz_2)$. As $\pol (\alpha )$ has  an
even-dimensional cell decomposition, $H^2(\pol(\alpha );\bbz)$ is free abelian
with rank equal of the dimension  of $H^2(\pol(\alpha );\bbz_2)$.
Therefore, the classes $c_i$ do not generate
$H^2(\pol(\alpha );\bbz )$ unless $\cals_m =
\{\emptyset\}$. \smallskip
 \decale{b)} In the proof of Proposition \ref{50} 
 the formula $c_m=-R$ could have been
 obtained directly from equation  (\ref{eqcm})
 and Proposition \ref{52}. The advantage of the previous argument is to be
 applicable to planar polygon spaces (see Proposition \ref{Stiefel}).
 \smallskip
 \decale{c)} By \cite[p. 109]{Fu}, the total Chern class of the
 (tangent bundle of the) upper path space is given by
$$c(\upol (\alpha ))\ = \ (1+R) \prod_{i=1}^m(1+U_i)\prod_{j=1}^{m-1}(1+V_j).$$
Using the relation $U_i=V_i+R$ and $(1+R)^2c(\pol (\alpha ))=c(\upol(\alpha))$
 gives
 $$c_1(\pol (\alpha))\  =\ (m-2)R + 2\sum_{i\in \cals_m}V_i\ 
 =\ \sum_{i=1}^m c_i .$$
 \smallskip
 \decale{d)} Using \ref{52} and \ref{51}, one gets the nicest expression for the
 cohomology class $[\omega ] \in H^2(\pol (\alpha ));\bbr ) $ of the symplectic
 form:
 $$ [\omega ]\  =\  \sum_{i=1}^m \alpha _i\,\, c_i.$$
 This is no surprise, since it is essentially how we calculated the $c_i$ in
 \cit{50}.
 \end{Remarks}
%
%
%

The great advantage of the $\{c_i\}$ over the generators $\{R, V_i\}$
is that they are
manifestly natural under permutation of the edges. Given a permutation
$\pi \in {\rm Sym}_m$, there is an isomorphism of $\pol(\alpha)$
and $\pol(\pi\alpha)$ given by reordering the steps. (This is a little
confusing in the polygon description, since one naturally thinks of
keeping
adjacent edges adjacent; instead one should simply think of a list of
$m$ vectors whose sum is zero, modulo rotation.)
>From the geometric construction of the $c_i$, one sees that under
$$ \pol(\alpha) \to \pol(\pi\alpha) $$
giving
$$ H^2(\pol(\pi\alpha)) \to H^2(\pol(\alpha)) $$
we have
$$ c_i \mapsto c_{\pi(i)}.$$

\begin{Proposition} \label{cgens}
The $\{c_i\}$ and Pontrjagin class $p$ are the generators in a manifestly
$S_m$-invariant presentation of the cohomology ring with coefficients in
$\bbz[1/2]$:
$$ H^*(\pol(\alpha); \bbz[1/2]) = \bbz[1/2][c_1,\ldots,c_m,p] / {\cal I}_c $$
 where $c_i$ is of degree 2 and $p$ of degree 4 and
  ${\cal I}_c$ is the ideal generated by the two families\smallskip
  \begin{tabbing} \renewcommand{\arraystretch}{5}
 \kern .7 truecm \= \+ (R1)\quad   \=
 $R^2\sum_{S\subset L \atop S\in \cals_m} \big(\prod_{i\in S}V_i )
 R^{|L-S|-1}$  \qquad \qquad \=    \= \kill

(R1)\> $c_i^2 - p$ \>  \> $i=1,\dots ,m-1$ \\ \bigskip

(R2)  \> $\displaystyle
   \sum^{\ }_{M\subset L \atop |M| \not\equiv |L| \bmod 2}
\big(\prod_{i\in M}V_i\big) p^{(|L-M|-1)/2}$  \>\> $L \in \call$
\end{tabbing}
\end{Proposition}

\preu
It is easiest to see this by returning to the original presentation in
\ref{40}. There were two steps necessary in Theorem \ref{41} to turn
this presentation
for the upper path space into one for the polygon space; (1) for each
$L \subseteq \{1,\ldots,m-1\}, L \in \call_m$ giving a relator (d),
subtract the corresponding relator (c) associated to $L \cup \{m\}$,
then (2) divide the difference by $R^2$ (now $c_m^2$).

Rewritten in terms of $c_i = R + 2 V_i = 2 U_i - R$ and $p$,
safely ignoring factors of 2, and performing the above two steps on (d),
 the relations (b)-(d) become
\begin{tabbing}
 \kern .7 truecm \= \+ (R1)\quad   \=
 $R^2\sum_{S\subset L \atop S\in \cals_m} \big(\prod_{i\in S}V_i )
 R^{|L-S|-1}$  \=    \= \kill

(b)  \> $(c_i+c_m)(c_i-c_m)$        \>\> $i = 1,\dots ,m-1$  \\

(c) \> $\prod_{i\in L} (c_i+c_m)$  \>\> $L\subset \{1,\dots ,m-1\}$ and $L\in \call_m$
\\ \smallskip
 (d')  \> $c_m^{-1}\,\big(\prod_{i\in L} (c_i-c_m)
                        -\prod_{i\in L} (c_i+c_m) \big)$ \\
   \>\>\> $L\subset \{1,\dots ,m-1\}$ and
$L\in \call$
\end{tabbing}

The relations (b) say that all the $c_i^2$ are equal,
which we knew in our ring,
since that's the Pontrjagin class $p$. Expand (c), pulling out
factors of $p$ where possible:
$$ \prod_{i\in L} (c_i+c_m) = \sum_{M \leq L} c_m^{|L-M|} \prod_{i\in M} c_i $$
$$
= \sum_{M \leq L \atop |M| \equiv |L| } p^{|L-M|/2} \prod_{i\in M} c_i
+ \sum_{M \leq L \atop |M| \not\equiv |L| } p^{(|L-M|-1)/2}
\prod_{i\in M \union \{m\}} c_i $$
(where the congruences are mod 2). Note that the products in both terms
of this last expression are over subsets of $L \union \{m\}$ with
odd complement. And in fact, every such subset appears this way exactly once;
this is relation (R2) in the case that the long subset of $\{1,\ldots,m\}$
contains $m$.

A similar analysis of (d') gives the relations (R2) for the case that
the long subset does not contain $m$; the negative terms in the first
expression cancel those subsets with even complement. \cqfd

This presentation is of most use in the case that $\pi \alpha = \alpha$,
and the induced isomorphism of polygon spaces is an automorphism.
In this case $\pi$ preserves the collection of long subsets,
and so permutes the relations given.
We emphasize that the generator $p$ is not necessary; its virtue is in
giving a much more efficient presentation.

\vskip 0.5 truecm
\section{Equilateral polygon spaces} \label{equil}

In this section we study the equilateral case, i.e. $\alpha_i =1$ for all $i$.
As usual, we require $\alpha$ to be generic, which in this case means exactly
that $m$ is odd. For the rest of this section we will use the notation
$\pol_m$ for the equilateral case with $m$ sides.

This space carries an action of $S_m$.
It is the one most commonly studied in algebraic geometry, because
the quotient
$\pol_m/S_m$ is a compactification of the moduli
space of $m$ {\em unordered} points in $CP^1$ -- in turn, the moduli space
of $m$-times-punctured genus zero algebraic curves.
Computing the cohomology ring of
this space is a classical problem, first solved by Brion \cite{Br}.

Since this space is an orbifold, it is most natural to consider its
rational\footnote{%
In fact one only has to invert primes up to $(m-1)/2$. The maximal
stabilizers of $\Sym_m$ come on triangles, and the very largest one
is $\Sym_{(m-1)/2} \times \Sym_{(m-1)/2} \times \Sym_{1}$.}
cohomology, particularly since one has a way to compute it:
$$ H^*(\pol_m/\Sym_m; \bbq) \iso H^*(\pol_m; \bbq)^{\Sym_m} $$
This requires one to understand the action of $\Sym_m$ on $H^*(\pol_m; \bbq)$
(first computed in \cite{Kl});
for our purposes, since we know it is generated in degree 2, we need only
understand the action on $H^2$. And this is easy, since the $\{c_i\}$
provide a basis demonstrating that the representation is the usual one
of $\Sym_m$ on $\bbq^m$ by permutation matrices.

In this section we show that the action of $\Sym_m$ on the {\em integral}
cohomology group $H^2(\pol_m)$ is {\em not} the standard one. And while
the set $\{c_i\}$ shows that it suffices to invert the prime $2$, it is
not necessary; in particular it suffices to invert the primes dividing $m$
(which, recall, is odd). This will find application in section \ref{plpol},
Theorem \ref{stdplanar}.

We then finish the section
by calculating the rational cohomology ring of $\pol_m/\Sym_m$,
by a different method than those in \cite{Br}\cite{KiPoly}.

Let $n$ be the product of the primes we are willing to invert.
To standardize the action on $H^2(\pol_m; \bbz[1/n])$,
we need to find a set $\{b_1,\ldots,b_m\} \in H^2(\pol_m; \bbz)$ satisfying
two criteria:

(C1) $\pi^*(b_i) = b_{\pi(i)}$ for all $i\in 1,\ldots,m$ and $\pi \in \Sym_m$

(C2) $\{b_1,\ldots,b_m\}$ is a basis of $H^2(\pol_m; \bbz[1/n])$.

(It suffices to take the $b_i$ in $H^2(\pol_m; \bbz)$, since one can
simply multiply them all by $n$ until this is true.)
To determine when $\{b_i\}$ is a basis, we make a matrix converting
from the known basis $\{R,V_i\}$ to $\{b_i\}$ and see if its determinant is
a unit in $\bbz[1/n]$.

We first find all solutions to criterion (C1). We look first for a vector
$b_m$ invariant under permutations of the first $m-1$ edges, and then
construct the other $m-1$ vectors from it
by applying the cycle $(123\ldots m)$.
Rationally it is easy to find all vectors invariant under permutations
of the first $m-1$ edges; take as the rational basis
$c_m$ and $\sum_{i=1}^{m} c_i$,
where $\{c_i\}$ are the Chern classes from section \ref{47}.

In our basis these are $-R$ and $(m-2) R + 2 \sum_{i=1}^{m-1} V_i$.
These two vectors do not span the intersection of the $\bbq$-space
they generate with the lattice $H^2(\pol_m; \bbz)$. To get all those vectors,
we must take combinations of the form
$$ b_m := \frac{x}{2} R + \frac{y}{2} ((m-2)R + 2 \sum_{i=1}^{m-1} V_i) $$
where $x \equiv y \bmod 2$ so that $R$ has an integer coefficient.

This basis was chosen to have easy transformation properties
under $(123\ldots m)$. The first cycles through the other $c_i = R + 2 V_i$,
whereas the second is fixed. For convenience of notation set $Z := m/2-1$.
The $\{b_i\} \to \{R,V_i\}$ transformation matrix
is then
$$ x/2 \pmatrix{
-1&&&&\cr
1&2&&&\cr
1&&2&&\cr
1&&&2&\cr
\vdots&&&&\ddots\cr}
+ y
\pmatrix{
Z&1&1&\cdots\cr
Z&1&1&\cdots\cr
Z&1&1&\cdots\cr
Z&1&1&\cdots\cr
&\vdots&&\cr
}
$$
We now compute the determinant.
Subtract rows $2$ through $m$ from their previous row, from the top down
(nothing propagates).
$$ x/2 \pmatrix{
-2&-2&&\cr
&2&-2&&\cr
&&2&-2&\cr
&&&\ddots\cr
1&&&&2\cr
}
+ y
\pmatrix{
&\vdots\cr
&0\cr
&0\cr
&0\cr
Z&1&1&\cdots\cr
}
$$
Subtract the first column from the second.
$$ x/2 \pmatrix{
-2&&&\cr
&2&-2&&\cr
&&2&-2&\cr
&&&\ddots\cr
1&-1&&&2\cr
}
+ y
\pmatrix{
&\vdots\cr
&0\cr
&0\cr
&0\cr
Z&1-Z&1&1&\cdots\cr
}
$$
Now add the second to the third, the third to the fourth, and so on
(this time things propagate).
$$ x/2 \pmatrix{
-2&&&\cr
&2&&&\cr
&&2&&\cr
&&&\ddots\cr
1&-1&-1&\cdots&1\cr
}
+ y
\pmatrix{
&\vdots\cr
&0\cr
&0\cr
&0\cr
Z&1-Z&2-Z&\cdots&(m-1-Z)\cr
}
.$$
This is a lower triangular matrix whose determinant is
$$ -x^{m-1} (x/2 + y(m-1 - (m/2-1))) = -x^{m-1} (x + ym)/2 $$
(recall that $m$ is odd, and $x \equiv y \bmod 2$, so this is actually
an integer).

\begin{Theorem} \label{caninvertm}
The action of $\Sym_m$ on the integral cohomology group $H^2(\pol_m)$
is never the standard permutation
representation. To standardize it it suffices to invert $2$ or $m$.
\end{Theorem}

\preu
At this point we are asking if $-x^{m-1} (x+ym)/2 $ can ever be a unit in
$\bbz$, that is to say, $\pm 1$. The first factor forces us to take $x=\pm 1$.
Then $y$ cannot be zero, since $x$ and $y$ have the same parity, so
$$ |ym+x| \geq |ym|-|x| = |y|m-1 \geq m-1 \geq 4 $$
so its quotient by $2$ cannot be as small as $\pm 1$. (Note that the
case $m=3$ fails for a simpler reason.)

To get the two possibilities advertised, take $x=2,y=0$ for $\bbz[1/2]$
(this is just the $c_i$ basis) and $x=m,y=1$ for $\bbz[1/m]$.
\cqfd

It is not too hard to find the exact conditions on $n$ making the
representation standard. The reader may find
it amusing to show that for $m=5$, it is necessary and sufficient that $n$
be divisible by a prime congruent to $0,2,3 \bmod 5$.

Hereafter in this section we work with rational coefficients,
and the $\{ c_i\}$ basis of $H^2(\pol_m; \bbq)$. In the equilateral case,
the presentation \cit{cgens} of the cohomology ring is particularly simple.
The minimal long subsets are exactly those with $(m+1)/2$ edges,
giving relators in degree $m-1$.

Focus first on the relators $(R1)\, c_i^2-p$, where $p$ was the extra
``generator'' in degree 4. These generate a sub-ideal agreeing with the
whole ideal up to degree $m-1$, and are easily seen to be a Gr\"obner
basis for this sub-ideal.

\begin{Lemma}
Let $\sigma_i$ be the $i$th symmetric polynomial in the $\{c_j\}$.
Then the subspace of $\bbq[p,\{c_i\}]$ generated by $p$ and the $\{\sigma_i\}$
maps onto the $\Sym_m$-invariant part of $H^*(\pol_m;\bbq)$.
For $*\leq m-3$ this map is an isomorphism.
\end{Lemma}

\preu
Consider the series of maps
$$
\big\{ \hbox{$\{ c_i^2-p\}$-reduced polynomials} \big\} \to
\bbq[p,\{c_i\}] \to
\bbq[p,\{c_i\}]/(c_i^2-p) \to
H^*(\pol_m; \bbq) .$$
The relations $\{c_i^2-p\}$ are easily seen to be
a Gr\"obner basis for the ideal they generate,
with respect to a reverse lexicographic order making powers of $p$ late
in the order.
Consequently the subspace of $\bbq[p,\{c_i\}]$ of reduced polynomials
maps isomorphically onto the quotient $\bbq[p,\{c_i\}]/(c_i^2-p)$ \cite{Ei}.
So the composition of the first two maps above is an isomorphism,
and therefore the composition of all three is an epimorphism.
In degrees below the omitted relations,
of degree $m-1$, the last map is an isomorphism.

The $\{c_i^2-p\}$-reduced polynomials are exactly combinations
of $p$ and the $\{c_i\}$ that are square-free in the $\{c_i\}$.
The condition of being square-free is preserved by the action of $\Sym_m$
permuting the $\{c_i\}$. So this composite map is
actually an $\Sym_m$-epimorphism,
and we can find the invariants up in our subspace rather than looking in
the quotient. These are exactly polynomials in $p$ and $\{c_i\}$ symmetric
in the $\{c_i\}$, which are generated by $p$ and the elementary symmetric
polynomials in $\{c_i\}$.
\cqfd

\begin{Corollary}
The $\Sym_m$-invariant part of $H^*(\pol_m;\bbq)$ is
generated by $\sigma_1$ of degree 2 and $p$ of degree 4,
with no relations up to degree $m-3$.
\end{Corollary}

\preu
The $\sigma_i$ are generated by $\sigma_1$ and $p$:
$$ \sigma_1 \sigma_i = (i+1) \sigma_{i+1} + (m-(i-1)) p \sigma_{i-1} $$
(where $\sigma_0 := 1$, $\sigma_{-1} = 0$). To see this, imagine multiplying a
product $\Pi_{i\in S} c_i$ by $c_j$. Either $j \notin S$, in which case the
product becomes one longer, or $j \in S$, in which case two $c_j$'s cancel
to become a $p$. The coefficients arise this way: in a product of $i+1$
things, any of them may be the new one, whereas in a product of $i-1$,
any one of the missing ones may be the one that just cancelled.
\cqfd

Since we now know the Betti numbers up to the middle dimension $m-3$,
by Poincar\'e duality we know all of them, and as in \cite{Br}
the Poincar\'e polynomial is quickly determined to be
$$ P_{\pol_m/\Sym_m} = \frac{(1 - t^{m-1}) (1 - t^{m+1})}{(1-t^2)(1-t^4)}. $$
In particular, since we know that there are only two generators (in
degrees 2 and 4), we know there are only two relations (in degrees
$m-1$ and $m+1$).

The relation in degree $m-1$ is not unexpected; it is the symmetric
combination of all the relations of degree $m-1$ in $H^*(\pol_m)$.
To find the one of degree $m+1$, we form for each $i=1,\ldots,m$ and $L \ni i$,
$|L|=(m+1)/2$ the S-polynomial of the corresponding relators in $(R1),(R2)$.
Our relation is the symmetric combination of those (one must check that
it is not $\sigma_1$ times the previous relator).
We omit the computations as the result is known \cite{Br}.

%

\begin{Theorem}
The rational cohomology ring of the equilateral polygon space
mod permutations is
$$ H^*(\pol_m/\Sym_m; \bbq) = \bbq[p,\{\sigma_i\}] / {\cal I} $$
where $p$ is of degree $4$, $\sigma_i$ of degree $2i$ for $i=0,\ldots,(m-1)/2$
and $\cal I$ is generated by the family
$$
(R1)\quad \sigma_1 \sigma_i = (i+1) \sigma_{i+1} + (m-(i-1)) p \sigma_{i-1}
 \qquad i=0,\dots,(m-3)/2
$$
and the two relators
$$\sum_{i \equiv (m-1)/2 \bmod 2}
        {m-i \choose (m+1)/2-i} p^{(\frac{m-1}{2}-i)/2} \sigma_i $$
and
$$\sum_{i \equiv (m+1)/2}
{\frac{m+1}{2} \choose i} / {m\choose i}
p^{\frac{1}{2}(\frac{m+1}{2}-i)} \sigma_i .$$
\end{Theorem}

It is worth explaining here exactly what problem Brion addresses, since it
is not obviously the one above. In both cases one is studying the action
of $SO(3) \times \Sym_m $ on $\prod_{i=1}^m S^2_{1}$. In the approach
above, we first perform the symplectic reduction by $SO(3)$, producing
the equilateral polygon space, and then take the quotient by $\Sym_m$.

One can perform these tasks in the opposite order. Regarding the $S^2$'s
as $CP^1$'s, one has available the celebrated homeomorphism\footnote{%
The orbifold structure is different, but this is not relevant for the
rational cohomology.}
of $(\prod_{i=1}^m CP^1) / \Sym_m$ with $CP^m$, taking the $m$ numbers to
their elementary symmetric combinations. This latter space is in turn the
projectivization of the $m+1$-dimensional irreducible representation
of $SU(2)$ -- on the projective space, the action factors through $SO(3)$.
It is this question -- the cohomology of $CP^m // SO(3)$ --
that Brion addresses and solves.

\vskip 0.5 truecm
\section{Planar polygon spaces}   \label{plpol}

 In this section, we study the {\it planar polygon space}:
 $$\polr (\alpha ) := \bigg\{(\rho _1,\dots,\rho _m)\in (\bbr^2)^m \biggm |
  |\rho _i|=\alpha_i
 \hbox{ and } \sum_{i=1}^m\rho _i=0  \bigg\}\biggm/O_2$$
 where $O_2$ acts on $(\bbr^2)^m$ diagonally. The more classical quotient
  by $SO_2$, denoted by
 $\pol (\alpha ;\bbr^2 )$, will also be considered.
 We assume  $\alpha $  generic, so the actions are free. The space
 $\polr (\alpha )$ is then a smooth manifold
 of dimension $m-3$ and $\pol (\alpha ;\bbr^2) \to \polr (\alpha )$ is
 a 2-fold cover.

 The $O(2)$-quotient $\polr (\alpha )$ is more natural for us because it is a
 submanifold of $\pol (\alpha )$. It can be interpreted  as a ``real part" of
 the K\"ahler manifold $\pol (\alpha )$: it is the fixed point set of the
 antiholomorphic involution $\rho\mapsto r\pcirc \rho$ where $r$ is the
 reflection $r(x,y,z)=(x,-y,z)$. More about that
  is to be found in  \cite[\SS 3 and 4]{HK}.
 The planar upper  path space $\upolr (\alpha )$ and abelian polygon space
 $\apolr (\alpha )$ are defined accordingly and can be seen as real parts of
 $\upol (\alpha )$ and $\apol (\alpha )$.

 We shall prove that a well known phenomenon for Grassmannians,
 toric manifolds, etc., also holds true for polygon spaces:

 \begin{Theorem} \label{cohppth}   Let $P$ (respectively: $P_\bbr$) stand for $\upol (\alpha )$,
 $\apol (\alpha )$ or $\pol (\alpha )$
  (respectively:  $\upolr (\alpha )$,
 $\apolr (\alpha )$ or $\polr (\alpha )$). Then there is
 is a ring isomorphism
 $$H^{2*}(P ; \bbz_2) \fl{\simeq}{} H^*(P_\bbr ;\bbz_2)$$
 sending elements of degree $2d$ to elements of degree $d$.
 \end{Theorem}


  For instance, for $P=\pol_r (\alpha )$, one gets:

  \begin{Corollary} \label{cohppcor}
   The cohomology ring of the planar polygon space $\polr(\alpha)$ with $\bbz_2$ as coefficient is
   $$H^*(\polr (\alpha );\bbz_2) = Z_2[R,V_1,\dots ,V_{m-1}] / {\cal I}_{\pol} $$
  where $R$ and $V_i$ are of degree 1
  and ${\cal I}_{\pol}$ is generated by the three families\smallskip

  \begin{tabbing} \renewcommand{\arraystretch}{5}
  \kern .7 truecm \= \+ (R1)\quad   \=
  $R^2\sum_{S\subset L \atop S\in \cals_m} \big(\prod_{i\in S}V_i\big)
  R^{|L-S|-1}$  \=    \= \kill

  (R1)\> $V_i^2+ RV_i$ \>  \> $i=1,\dots ,m-1$ \\ \bigskip

  (R2)  \> ${\displaystyle \prod^{\ }_{i\in L} V_i}$  \>\>
  $L\in \call_m$ \\ \bigskip

  (R3)  \> ${\displaystyle \sum^{\ }_{S\subset L \atop S\in \cals_m} \big(\prod_{i\in
  S}V_i\big) R^{|L-S|-1}}$  \>\> $L\subset \{1,\dots ,m-1\}$ and $L\in \call$
  \end{tabbing}
  \end{Corollary}

{\sc Proof of \ref{cohppth}: }
 As seen in Section \ref{definitions}, the manifolds
 $\upol (\alpha )$ and $\apol (\alpha )$ are toric manifolds.
 Therefore, Theorem \ref{cohppth} is true by \cite[Theorem 4.14]{DJ}
 and a proof is only required for $\pol (\alpha )$.

 We first establish that for each $k\in\bbn$: \setcounter{equation}{0}
 \begin{equation}
 {\rm dim\,} H^k (\polr (\alpha );\bbz_2) \leq
 {\rm dim\,} H^{2k} (\pol (\alpha );\bbz_2)\label{ineq}
 \end{equation}
 where ${\rm dim\,} $ means the dimension as a vector space over the field $\bbz_2$.
 This is done by induction on the number $m$ of edges. The statement is trivial
 for $m=3$ where   $\polr (\alpha )=\pol (\alpha ) = \hbox{one point}$. It is
 also obviously true for $m=4,5$ by the list of all polygon spaces (see \cite[Section
 6]{HK}).

 We use the notations of the proof of Lemma \ref{polevcoh}.
 By \cite[Theorem 3.2]{Ha}, the diagonal-length function $\delta: \pol (\alpha
 )\fl{}{}
 \bbr $ given by   $\delta (\rho ):=|\rho  _m - \rho  _{m-1}|$ is a Morse-Bott
 function on $\polr (\alpha )$. The critical
 points are the same as those for $\pol (\alpha )$ but, for each of them,  the index is divided by
 2. They are isolated except possibly for the  two extrema.
 The pre-image $M_{\rm max}$ of the maximum is either a point or $\polr(\alpha
 _1,\dots ,\alpha _{m-2},\alpha _m +\alpha _{m-1})$. For the pre-image
 $M_{\rm min}$ of the
 minimum, there are three possibilities:
 \decale{--} one point
 \decale{--} $\polr(\alpha _1,\dots ,\alpha _{m-2},\alpha _m
 -\alpha _{m-1})$
 \decale{--} a circle bundle over
 $\polr(\alpha _1,\dots ,\alpha _{m-2},\alpha _m -\alpha _{m-1})$ (when the
 minimum is $0$).

 By induction on $m$, inequality (\ref{ineq}) holds for
 $M_{\rm min}$, $M_{\rm max}$ and $\polr (\alpha ) - M_{\rm max}$.
 As in Proposition \ref{ex-seq} one one gets an exact sequence
 \begin{eqnarray} \nonumber
 \fl{}{}H^{*-1}(M_{\rm max}) \fl{}{} H_{n-*}(\polr(\alpha )-M_{\rm max})
  \fl{}{}\\ \fl{}{}  H^*(\polr (\alpha ))\fl{i^*}{}  H^*(M_{\rm max})\fl{}{}
 \label{seq3}
 \end{eqnarray}
  For $\pol (\alpha )$ this exact sequence is cut into short ones by
  Proposition \ref{ex-seq}. This enables us to propagate  inequality
  (\ref{ineq}) to $\polr (\alpha )$.

 As in the proof of Theorem \ref{64}, the class
 $R^2\in H^*(\upolr (\alpha );\bbz_2)$ is Poincar\'e dual
 to $\polr (\alpha )$. By the proof of Proposition \ref{ker}, the
 annihilator $\Ann (R^2)$ in $H^*(\upolr (\alpha );\bbz_2)$ of the cup product
 with  $R^2$
 contains the kernel $\ker i^*$ of
 $i^*: H^*(\upolr (\alpha );\bbz_2)\to H^*(\polr (\alpha );\bbz_2)$. Combining
 with inequality (\ref{ineq}) gives the following sequence of inequalities:

 \begin{eqnarray*}
 \dim \big(H^*(\upolr (\alpha );\bbz_2)\big /\Ann (R^2)\big) & \leq & \\
  \dim \big(H^*(\upolr (\alpha );\bbz_2)\big /\ker i^* \big) &  = & \dim {\rm
  Image\,}i^* \leq \\
   \dim H^*(\polr (\alpha );\bbz_2) & \leq \\
    \dim H^*(\pol (\alpha );\bbz_2) & =  &
    \dim \big(H^*(\upolr (\alpha );\bbz_2)\big/\Ann (R^2)\big) ,
     \end{eqnarray*}
 the last equation being Theorem \ref{62}.  The two ends being equal, all the
 above inequalities are equalities. Therefore $\Ann (R^2)=\ker i^*$,
 $i^*$ is surjective,
 ${\rm dim\,} H^k (\polr (\alpha );\bbz_2) =
 {\rm dim\,} H^{2k} (\pol (\alpha );\bbz_2) $ and one has an isomorphism
 $$H^* (\polr (\alpha );\bbz_2) \simeq H^*(\upolr (\alpha );\bbz_2))\big /\Ann
 (R^2).$$
 This proves Theorem \ref{cohppth}.   \cqfd

 We now turn our attention to the 2-fold cover
 $\kappa : \pol (\alpha ;\bbr^2 )\to \polr(\alpha)$. Seen as a principal  $O_1$ cover, it is
 determined by its Stiefel-Whitney class
 $w_1(\kappa)\in H^1(\polr(\alpha);\bbz_2)$.

 \begin{Proposition} \label{Stiefel} $w_1(\kappa) = R$.
 \end{Proposition}

 \preu  The $O_1$-bundle $\kappa$ is the planar analogue of $U_1$-bundle $\xi_m$ introduced
 in Section \ref{47}. The proof that $c_1(\xi_m)=-R$ (proof of Proposition
 \ref{50}) then works $\bmod 2$ to give $w_1(\kappa) = R$. \cqfd

Lastly, we discuss equilateral planar polygons. We cannot say much about
the quotient by the symmetric group since those calculations involve
inverting the prime 2. There is something to say about the action:

\begin{Proposition} \label{stdplanar}
The action of $S_m$ on $H^2(\polr_m; \bbz_2)$
is the standard one on $\bbz_2^m$.
\end{Proposition}

\preu
Take the $\bbz[1/m]$-basis of $H^2(\pol_m)$ from Theorem \ref{caninvertm}.
This becomes a basis once $m$ is invertible, which it is over $\bbz_2$.
\cqfd

\vskip 0.5 truecm
\section{Examples} \label{Exples}

 \begin{ccote} \label{exa} \rm  Suppose $\cals_m (\alpha )=\{\emptyset\}$, for example if
 $\alpha = (1,\dots ,1,m-2)$.  It follows from \cite[Proposition (4.1)]{Ha} that $\pol (\alpha )$ is diffeomorphic
 to the complex projective space
 $\bbc P^{m-3}$. Knowing this, we can test our different results for the homology
 or cohomology of $\pol (\alpha )$.

 As $\cals_m (\alpha )=\{\emptyset\}$, the expression of the Poincar\'e
 polynomial $P_{\pol (\alpha )}$ given in Theorem \ref{38} is a 1-term sum:
 $$P_{\pol (\alpha )} = \frac{1-t^{2(m-2)}}{1-t^2} = 1 + t^2 + \cdots +
 t^{2(m-3)}$$
 which is indeed the Poincar\'e polynomial of $\bbc P^{m-3}$. Observe that the
 formula for $P_{\pol (\alpha )}$ in terms of $\cals$  given in
 Remark \ref{38b} would have $2^{m-1}-1$ terms!

 For the cohomology ring $H^*(\pol (\alpha )))$, Theorem \ref{62} asserts that it
 is the quotient of $\bbz[R,V_1,\dots , V_{m-1}]$ by an ideal $\cali$ generated
  by the families
 of relators (R1), (R2) and (R3). As $\{i\}\in\call_m$ for all $1\leq i\leq m-1$,
 all the $V_i$'s are killed by (R2) and (R1) becomes empty.
  Family (R3) contains one
 element, for $L = \{1,\dots ,m-1\}$. As $\cals_m (\alpha )=\{\emptyset\}$,
 this relator is $R^{m-2}$. Thus $H^*(\pol (\alpha )) = \bbz [R]/(R^{m-2})$, the
 cohomology ring of $\bbc P^{m-3}$.

 In the planar case, one has $\polr (\alpha ) \simeq \bbr P^{m-3}$ and
 $\pol (\alpha ;\bbr^2) \simeq S^{m-3}$.
  \end{ccote}

 \begin{ccote}\rm\label{exaa}  Consider the case where
 $\cals_m (\alpha )$ contains $ \{1,\dots ,m-3\}$
(for instance: $\alpha =(\varepsilon ,\dots ,
 \varepsilon , 1,1,1)$ with $(m-3)\varepsilon<1$). Then, any $r\in \pol (\alpha )$ has a unique
representative $\rho$ with $\rho_m = (1,0,0)$ and $\rho(m-1) =
 (x,y,0)$ with $y>0$. The class $r$ is then determined by $\rho(1),\dots
 ,\rho(m-3)$ and there is no constraint on these vectors. Therefore, $\pol (\alpha
 )$ is symplectomorphic to $\prod_{i=1}^{m-3} S^2_{\alpha _i}$. In the planar
 case,  $\polr (\alpha )$ is diffeomorphic to $\prod_{i=1}^{m-3} S^1$. The space
 $\pol (\alpha ;\bbr^2$ is not connected:
 $\pol (\alpha ;\bbr^2) \simeq S^0 \times \prod_{i=1}^{m-3} S^1$.

 Let us compute the cohomology ring $H^*(\pol (\alpha )) =
 \bbz[R,V_1,\dots ,V_{m-1}]/\cali$.
 One has $\cals_m = \Delta^{m-4}$ and the minimal elements of $\call_m$ are the
 singletons $\{m-2\}$ and  $\{m-1\}$. Therefore, relators (R2) reduce to
 $V_{m-2}$ and $V_{m-1}$. The minimal $L\subset \{1,\dots ,m-1\}$
 in $\call$ is $L=\{m-2,m-1\}$. For this $L$, relator (R3) is $R$.
 The other relators of the family (R3) all have $R$ as a factor. Finally,
 using (R1), one finds
 $$ H^*(\pol (\alpha )) =  \bbz[V_1,\dots ,V_{m-3}]\big / (V_1^2, \dots ,V_{m-3}^2)$$
 as expected.

In the particular case $m=3$, the cohomology ring reduces to the degree 0 part
(no wonder since a triangle space is just a point).
\end{ccote}

 \begin{ccote}\rm\label{quadr} Consider the two cases of quadrilaterals mentioned in Section \ref{shortlong}:
  $\alpha = (1,1,1,2)$ and $\alpha' = (1,2,2,2)$. As $\cals_4(\alpha)=\{\emptyset\}$, we are in case
 \ref{exa} and $H^*(\pol (\alpha )) = \bbz [R] / (R^2)$.  The case $\alpha '$ is like example \ref{exaa} and
 $H^*(\pol (\alpha' )) = \bbz [V_1] / (V_1^2)$ (in particular, $R=0$). Therefore,
 $\xi(\alpha )$ is the non-trivial $SO_3$-bundle over $S^2$ whereas
 $\xi(\alpha ')$ is the trivial one. In the same way,
 $\pol (\alpha ;\bbr^2)\to \polr (\alpha )$
 is the connected 2-fold cover of $S^1$ whereas $\pol (\alpha ';\bbr^2)\to \polr (\alpha' )$
 is the trivial cover.
\end{ccote}

 \begin{ccote}\rm\label{exb} {\it The regular pentagon: } $\alpha =(1,1,1,1,1)$. The space
 $\pol (\alpha )$ is a smooth manifold diffeomorphic to
 $(S^2\times S^2)\sharp 3 \overline {\bbc P}^2 \simeq
    {\bbc P}^2\sharp 4 \overline {\bbc P}^2$ \cite{Kl}
\cite[(6.3)]{HK}.

 One has
 $\cals_5=  \{\{1\},\{2\},\{3\},\{4\}\}$ and therefore $H^*(\pol (\alpha )))$ is
 generated by $R, V_1,\dots ,V_4 \in  H^2(\pol (\alpha ))$.
 The minimal elements of $\call_5$ are
 the doubletons $\{i,j\}$ for $i,j = 1,2,3,4$; hence family (R2) is
 generated by relators $V_i V_j$.
 The subsets of $\{1,2,3,4\}$
 which are elements of $\call $ are $L_j:=\{1,2,3,4\}-\{j\}$ and
 $L:=\{1,2,3,4\}$. This gives rise to five relators in family (R3):
 $$\begin{array}{lll}
 (1)\quad  L=\{1,2,3\} & : \kern .6 truecm & R^2 + RV_1 + RV_2 + RV_3 \\
 (2)\quad L=\{1,2,4\} & :  & R^2 + RV_1 + RV_2 + RV_4 \\
 (3)\quad L=\{1,3,4\} & :  & R^2 + RV_1 + RV_3 + RV_4 \\
 (4)\quad L=\{2,3,4\} & :  & R^2 + RV_2 + RV_3 + RV_4 \\
 (5)\quad L=\{1,2,3,4\} & :  & R^3 + R^2V_1 + R^2V_2 + R^2V_3 + R^2V_4 \ .
  \end{array}$$
  One deduces that $RV_1 = RV_2 =RV_3=RV_4$, and
  $R^2=-3RV_1$. One also has relators (R1) : $V_i^2+RV_i$. One then checks that
  everything in degree 3 vanishes. Let us take $T=R+V_1+V_2+V_3+V_4$ and the
  $V_i$'s as a basis for $H^2(\pol (\alpha ))$ and $RV_1$ as a basis of
  $H^4(\pol (\alpha ))$. With these bases, the cup product
  $H^2(\pol (\alpha ))\times H^2(\pol (\alpha )) \to H^4(\pol (\alpha ))$ is given
  by the following matrix.
  $$\pmatrix{1 & 0 & 0 & 0 & 0 \cr
             0 & -1 & 0 & 0 & 0 \cr
             0 & 0 & -1 & 0 & 0    \cr
             0 & 0 & 0 & -1 & 0       \cr
             0 & 0 & 0 & 0 & -1 \cr}$$
  which is indeed the intersection form of  ${\bbc P}^2\sharp 4 \overline {\bbc
  P}^2$.

  By Proposition \ref{52}, the class $[\omega]\in H^2(\pol (\alpha ))$ of the symplectic form $\omega$ is
 $$ [\omega] = 3R + 2V_1 + 2V_2 + 2V_3 + 2V_4 $$ and therefore $ [\omega]^2 = 5
 RV_1$. The Liouville volume $\int_{\pol (\alpha )}\omega^2/2$ is then $5/2$.
 We  get exactly the area of the ``moment polytope" of \cite[Figure 3]{HK}
 (we put ``moment polytope" between quotes since the regular
 pentagon space is only a
 limit case of toric manifold; see \cite[(6.3)]{HK}). This illustrates the
 Duistermaat-Heckmann theorem.
 \end{ccote}

 \begin{ccote}\rm\label{exbb}  The pentagon spaces for generic $\alpha $'s are all classified
 \cite[(6.2)]{HK}. They are toric manifolds and thus classified by their
 moment polytope. The reader can check, as for the regular pentagon space, that
 one gets the correct intersection forms for these 4-manifolds and that the
 Liouville volume is the area of the moment polytope.
 \end{ccote}

 \begin{ccote}\rm\label{exc} Consider the hexagon spaces $\pol (\alpha )$ and $\pol (\alpha ')$
 for
 $$ \alpha := (2,2,3,5,5,10) \ \hbox{ and } \alpha ' :=  (2,2,3,5,5,8).$$
 One has $\cals_6(\alpha )= \{\emptyset,\{1\},\{2\},\{3\}\}$ and
 $\cals_6(\alpha' )= \{\emptyset,\{1\},\{2\},\{3\},\{1,2\}\}$.
 Using \ref{38}, one sees that these polygon spaces cannot be distinguished by their Poincar\'e polynomial:
 $$P_{\pol (\alpha)}=P_{\pol (\alpha')}= 1 + 4 t^2 + 4t^4+ t^6 \ .  $$
 As for the cohomology rings, we deduce from Theorem \ref{62} that
 $H^*(\pol (\alpha ))$ has the following description

 $$\begin{array}{lll}
                     & \hbox{generators:} &  \hbox{further relations:} \\
  H^2(\pol (\alpha ))\quad &  R,V_1 , V_2 , V_3 \\
   H^4(\pol (\alpha ))) &  R^2, V_1^2, V_2^2 , V_3^2\ \quad &
   V_1V_2 = V_1V_3 =V_2V_3 = 0 \, , \, V_i^2=-RV_i \\
    H^6(\pol (\alpha )) &   V_1^3 = V_2^3= V_3^3 &
    2V_1^3 =-R^3\,
   \end{array} $$

   whereas, for $H^*(\pol (\alpha '))$, one has

 $$\begin{array}{lll}
                     & \hbox{generators:} &  \hbox{further relations:} \\
  H^2(\pol (\alpha '))\quad &  R,V_1 , V_2 , V_3 \\
   H^4(\pol (\alpha' )) &  R^2, V_1^2, V_2^2 , V_1V_2\ \quad &
    V_1V_3 =V_2V_3 = 0 \, , \, V_i^2=-RV_i\, , \, V_3^2=-R^2 \\
    H^6(\pol (\alpha ')) &  V_3R^2=R^3  & V_1^3 = V_2^3= 0 ,\, RV_1V_2 = - R^3
   \end{array} .$$

 These two rings are nonisomorphic, even over $\bbz_4$.
 It is a computer algebra exercise to show that
 in $H^*(\pol (\alpha ))\otimes \bbz_4$ there are 72 elements $x$ with
 $x^3=0$ whereas this number is 80 for $H^*(\pol (\alpha' ))\otimes \bbz_4$.
 One can check by hand the more relevant fact that there is no
 isomorphism between $H^* (\pol (\alpha ))$ and $H^* (\pol (\alpha ))$
 preserving the classes
 $R$'s. Indeed, $R^3 \bmod 2$ generates $H^6(\pol (\alpha ');\bbz _2)$ whereas
 $R^3=0$ in  $H^6(\pol (\alpha) ;\bbz _2)$.  But, by \ref{49} and \ref{50},
 $R \bmod 2$ is the second
 Stiefel-Whitney class of the $SO_3$-bundle $\xi (\alpha )$
 defined in Section \ref{47}. In particular, $A(\alpha )$ and $A(\alpha ')$ are not
 $SO_3$-equivariantly diffeomorphic.
 \end{ccote}


\vskip .5 truecm
\small

 \noindent Jean-Claude HAUSMANN\\ Math\'ematiques-Universit\'e\\ B.P. 240, \\
 CH-1211 Gen\`eve
 24, Suisse\\ hausmann@math.unige.ch
 \vskip 0.3 truecm

 \noindent Allen KNUTSON \\
 Department of Mathematics\\
 Brandeis University\\
 Waltham, MA 02254-9110 USA\\
 allenk@alumni.caltech.edu

\end{document}